\documentclass[a4paper,11pt]{article}

\usepackage{amsfonts}
\usepackage{amsmath}
\usepackage{amssymb}
\usepackage{graphicx}
\usepackage{amscd}
\usepackage{tikz}
\usepackage{epsfig}
\usepackage{psfrag}
\usepackage[margin=2cm]{geometry}

\newtheorem{theorem}{Theorem}

\newtheorem{conjecture}[theorem]{Conjecture}

\newtheorem{definition}{Definition}

\newtheorem{lemma}{\indent Lemma}

\newtheorem{proposition}{Proposition}

\newcommand{\qp}[2]{(#1;q)_{#2}}
\newcommand{\qpp}[1]{(q;q)_{#1}}

\def\max{\mathop{\rm max}}

\def\Z{\mathbb{Z}}
\def\N{\mathbb{N}}

\def\P{\mathcal{\widetilde{P}}}
\title{Generalized knots-quivers correspondence}
\date{}
\author{Marko Sto\v si\'c\footnote{
CAMGSD, Departamento de Matem\'atica, 
Instituto Superior T\'ecnico,
Av. Rovisco Pais 1, 1049-001 Lisbon, Portugal, and
Mathematical Institute SANU, Knez Mihajlova 36, 11000 Belgrade, Serbia. 
(e-mail: {\tt mstosic@isr.ist.utl.pt})
}}
\begin{document}
\maketitle

\abstract{We propose a generalized version of knots-quivers correspondence, where the quiver series variables  specialize to arbitrary powers of the knot HOMFLY-PT polynomial series variable. We explicitely compute quivers for large classes of knots, as well as many homologically thick 9- and 10-crossings knots, including the ones with the super-exponential growth property of colored HOMFLY-PT polyomials. In addition, we propose a new, compact, quiver-like form for the colored HOMFLY-PT polynomials, where the structure of colored differentials is manifest. In particular, this form partially explains the non-uniqueness of quivers corresponding to a given knot via knots-quivers correspondence.}


\section{Introduction}
Knots-quivers correspondence has been introduced in \cite{KRSS}. To each knot this correspondence associates a quiver so that the generating series of colored HOMFLY-PT polynomials of a knot matches the quiver generating series under suitable variable specializations. These variable specializations are necessary since HOMFLY-PT generating series has a single series variable, $x,$ whereas quiver generating series has as many distinct variables, $x_1,\ldots,x_m$, as the number of its nodes. Usually the variable specializations take the form
\begin{equation}\label{x}x_i=(-1)^{s_i}a^{a_i}q^{q_i}x,\quad i=1,\ldots,m,\end{equation} for suitable integers $a_i,$ $q_i$, $s_i$. Such specializations were used in the original paper \cite{KRSS}, as well as in \cite{SW1,SW2} to prove knot-quivers correspondence and find explicit quivers for large classes of knots. Moreover, other combinatorial interpretations of the coefficients of the classical limits of such quivers were obtained, together with the explicit closed formulas for them. \\

Notable feature of specializations (\ref{x}), is the fact that all specializations are linear, i.e. the exponent of $x$ in RHS is always equal to 1. In recent results, where the other, knot complement, invariants $F_K$ were expressed in a quiver series form (\cite{EGG+,EKL2}), it turned out that the specializations with the  higher powers  of $x$ were necessary. That is, the series of the invariant $F_K$ could be written in a quiver series form, but where the specializations of the  variables of the quiver generating series have the form
\begin{equation}\label{x1}x_i=(-1)^{s_i}a^{a_i}q^{q_i}x^{n_i},\quad i=1,\ldots,m,\end{equation}
where $n_i$ are nonnegative integers. In particular some of $n_i$ can be larger than $1$. The generating series obtained in this way have also been called {\it{generalized quiver partition functions}} in \cite{EKL}, where they have been related to knot homologies obtained by holomorphic curves count.

One of the reasons to extend the knots-quivers correspondents so that the quiver generating series can have higher level generators, is the super-exponential growth behaviour of the colored HOMFLY-PT polynomials. Namely, the number of terms of the $r$th colored HOMFLY-PT polynomial tend to grow faster than exponentially on $r$. Although many smaller knots did have exponential growth and consequently the linear specialization (\ref{x}) were enough, for larger knots we need the specializations (\ref{x1}). The goal of this paper is to extend the original knots-quivers correspondence from \cite{KRSS}, and {{this new correspondence}} we call {\it{generalized knots-quivers correspondence}}.

\begin{conjecture}[Generalized knots-quivers correspondence]
For a given knot the generating function of its appropriately normalized symmetrically colored HOMFLY-PT polynomials, can be written as the quiver generating series of a suitable symmetric quiver with the specialization of variables of the form (\ref{x1}).
\end{conjecture}

Remarkably, it turns out that the colored HOMFLY-PT homologies can be recovered from the quiver. In this setting, in various examples that we give in the paper,  the generators and the homological degrees of colored HOMFLY-PT homologies can be recovered from the quiver.
In particular, we present this in detail for the knot $\mathbf{9_{42}}$.

Another important thing is that many of the quivers that we will present are obtained by using ansatz, where we are looking first for expressions of a special quiver-like form, where the corresponding matrices are much smaller in size. Apart from making this computationally much simpler problem, this specific form incorporates the structure of the differentials and colored HOMFLY-PT homologies, and in particular, perfectly matches the conjectures of the structure of such colored HOMFLY-PT homologies \cite{GGS,GS}. Remarkably, we manage to find quivers of such a specific form for all the knots that we studied.

As in the case of quivers with only level one generators, here we also have that quivers corresponding to a given knot are non-unique. Moreover, when nodes of higher level are allowed, there are even more possibilities for quivers having the same generating series.  On the other hand, restricting only to quivers obtained in specific form that we observe in this paper, in many ways reduce the non-uniqueness, and explain the rich, discrete structure of the quivers with the minimal size, including (large parts of) the permutahedron structure studied in \cite{permutahedron}.\\

In Section \ref{not} we will give precise notation and statements of the generalized knots-quivers correspondence. We also give a particular form of the quiver that is suitable for our computations. In Section \ref{results} we list the explicit results for many knots. In Section \ref{sec942} we give more details about the flagship example $\mathbf{9_{42}}$. In Section \ref{homology} we give extensions and relationship with colored HOMFLY-PT homologies.\\

The results of this paper are collection of various results and computation from the past couple of years. I greatly benefited, both by discussions with numerous experts in the field, as well as by sharing their computations of different knot invariants. In particular, I am thankful to Ben Cooper, Tobias Ekholm, Sergei Gukov,  Piotr Kucharski, Satoshi Nawata, Piotr Su\l kowski, and Paul Wedrich. \\

\textbf{Acknowledgments:} The author was supported by FCT through project "Knots, quivers, and categorification", DOI number 10.54499/2020.02453.CEECIND/CP1587/CT0007. This work was also partially supported by the Science Fund of the Republic of Serbia, Grant GWORDS -- ``Graphical Languages", No. 7749891.

\section{Notation}\label{not}

Generating function of all symmetrically colored unreduced {HOMFLY-PT} polynomials of a knot $K$ is given by:
\begin{equation}\label{kser}
P(K)(x,a,q)=\sum_{r\ge 0} \bar{P}_r(K)(a,q) x^r,
\end{equation}
where $\bar{P}_r(K)(a,q)$ is the unreduced $Sym^r$-colored HOMFLY-PT polynomials of $K$.

For a quiver $Q$ with $m$ vertices, let $C$ be a matrix of a quiver. The corresponding quiver generating series 
is given by:
\begin{equation}
P_C(x_1,\ldots,x_m)=\sum_{d_1,\ldots,d_m} \frac{(-q^{\frac{1}{2}})^{\sum_{i,j=1}^m C_{i,j}d_id_j}}{(q;q)_{d_1}\cdots(q;q)_{d_m}} x_1^{d_1}\cdots x_m^{d_m}. \label{P-C}
\end{equation}

Knots-quivers correspondence was originally postulated in \cite{KRSS,KRSSlong}. For a given knot $K$ it conjectures the existence of a quiver $Q$ such that the corresponding quiver generating series (\ref{P-C}) matches the generating series (\ref{kser}) of the colored HOMFLY-PT of the knot $K$, after specializations:
\begin{equation}\label{xx}x_i=(-1)^{s_i}a^{a_i}q^{q_i}x,\quad i=1,\ldots,m,\end{equation} for suitable integers $a_i,$ $q_i$, and $s_i$, for $i=1,\ldots,m$.
This can be restated explicitely in the following way for the reduced $Sym^r$-colored HOMFLY-PT polynomials:

\begin{definition}[Quiver-form of colored HOMFLY-PT polynomials]\cite{KRSSlong}
For a knot $K$ we shall say that its reduced symmetrically colored HOMFLY-PT polynomials, $P_r(K)$ can be written in a quiver form if there is a positive integer $m$, a symmetric integral $m\times m$ matrix $C$, and two $m$-tuples of integers $\mathbf{a}=(a_1,\ldots,a_m)$ and $\mathbf{q}=(q_1,\ldots,q_m)$, such that for every $r\in \N$:
 \begin{eqnarray}
P_r(K)(a,q) = \sum_{{\tiny{\begin{array}{c}d_1,\ldots,d_m\ge 0\\
\sum_i {d_i=r}\\
\end{array}}}}  \frac{(q;q)_r}{\prod_{i=1}^m {(q;q)_{d_i}}}  (-1)^{\sum_i C_{ii} d_i} a^{\sum_i a_i d_i} q^{\sum_i q_i d_i} q^{\frac{1}{2}(\sum_{i,j} C_{ij}d_id_j-\sum_i C_{ii} d_i)}  
\label{kq1}
\end{eqnarray}
\end{definition}

The generalized version we propose is the following:
\begin{definition}[Generalized quiver-form of colored HOMFLY-PT polynomials]
For a knot $K$ we shall say that its reduced symmetrically colored HOMFLY-PT polynomials, $P_r(K)$ can be written in a quiver form if there is a positive integer $m$, a symmetric integral $m\times m$ matrix $C$, two $m$-tuples of integers $(a_1,\ldots,a_m)$ and $(q_1,\ldots,q_m)$, and an $m$-tuple of nonnegative integers $(n_1,\ldots,n_m)$, such that for every $r\in \N$:
 \begin{eqnarray}
P_r(K)(a,q) = \sum_{{\tiny{\begin{array}{c}d_1,\ldots,d_m\ge 0\\
\sum_i {n_i d_i=r}\\
\end{array}}}}  \frac{(q;q)_r}{\prod_{i=1}^m {(q;q)_{d_i}}}  (-1)^{\sum_i C_{ii} d_i} a^{\sum_i a_i d_i} q^{\sum_i q_i d_i} q^{\frac{1}{2}(\sum_{i,j} C_{ij}d_id_j-\sum_i C_{ii} d_i)}  
\label{kq2}
\end{eqnarray}
\end{definition}

This corresponds to the version of knots-quivers corresponds, where the specialization of the quivers series variables are such that the they are not necessarily linear in $x$, i.e., instead of (\ref{xx}) one has 
\begin{equation}\label{xx1}x_i=(-1)^{s_i}a^{a_i}q^{q_i}x^{n_i},\quad i=1,\ldots,m,\end{equation}
where $n_i$, $i=1,\ldots,m$ are nonnegative integers.\\

\subsection{LMOV conjecture}
One of the original goals of knots-quivers correspondence was to prove the LMOV conjecture of knots, and to compute and refine their BPS numbers, for more details see \cite{KRSSlong}. Briefly, the LMOV conjecture (for symmetric representations) can be stated in the following way: first one re-writes the colored HOMFLY-PT generating series (\ref{kser}) for a knot $K$ as

$$P(K)(x,a,q):=\sum_{r\ge 0} \bar{P}_r(K)(a,q) x^r =\exp \left(  \sum_{n,r\ge 1} \frac{1}{n} f_r(a^n,q^n) x^{rn}  \right),$$
where
$$ f_r(a,q)=\sum_{i,j} \frac{N_{r,i,j}a^iq^j}{q-q^{-1}}. $$
Then the LMOV conjecture states the numbers $N_{r,i,j}$ (called BPS numbers) are integers.\\

On the other hand, in the case of a symmetric quiver $Q$,  the motivic Donaldson-Thomas invariants  $\Omega_{d_1,\ldots,d_m;j}$ of $Q$, \cite{KS}, can be interpreted as the intersection Betti numbers of the moduli space of all semisimple representations of $Q$, or as the Chow-Betti numbers of the moduli space of all simple representations \cite{FR,MR}. These invariants are encoded in the following product decomposition of the quiver generating (\ref{P-C}) series of the quiver $Q$
\begin{equation*}
P_C(x_1,\ldots,x_m)=
\prod_{(d_1,\ldots,d_m)\neq 0} \prod_{j\in\mathbb{Z}} \prod_{k\geq 0} \Big(1 -  \big( x_1^{d_1}\cdots x_m^{d_m} \big) q^{\frac{j+1}{2}+k} \Big)^{(-1)^{j+1}\Omega_{d_1,\ldots,d_m;j}}.   \label{PQx-Omega}
\end{equation*}

\begin{theorem}\cite{Ef} \label{efimov}$\Omega_{d_1,\ldots,d_m;j}$ are nonnegative integers. \end{theorem}


In the same was as in the original paper \cite{KRSSlong}, for all the knots for which there exists a quiver satisfying generalized knots-quivers correspondence (\ref{kq2}), the LMOV conjecture holds.

\begin{proposition}
For all knots satisfying the generalized knots-quivers correspondence, the LMOV conjecture holds.
\end{proposition}

Indeed, the BPS numbers can again be written as integral linear combination of the motivic DT invariants of the corresponding quiver, and hence are integers by Theorem \ref{efimov}.

\section{Special form of colored HOMFLY-PT polynomials}\label{special}

In this section we present explicit results for the symmetric colored HOMFLY-PT polynomials of various knots, written in a quiver-like form.
To that end, we shall look for the quivers of a particular special form.

More precisely, we rewrite the generating series of the reduced $S^r$-colored HOMFLY-PT polynomials in the following way:
\begin{eqnarray}
P(K)(x;a,q)&=&\sum_{r} P_r(K)(a,q) x^r= \nonumber\\
&=&\sum_{d_1,\ldots,d_k} (-1)^{\sum_i C_{ii} d_i}q^{\frac{1}{2}(\sum_{i,j} C_{ij}d_id_j-\sum_i C_{ii} d_i)} (q;q)_{\sum_i n_i d_i}\prod_{i=1}^k \frac{x_i^{d_i}}{(q;q)_{d_i}} \times \nonumber\\
&&\quad\quad\quad\times(a^{-1}q;q^{-1})_{\sum_i r_i d_i} (a^{-1}q^{-\sum_i n_i d_i};q^{-1})_{\sum_i l_i d_i},\label{specform}
\end{eqnarray}
where $C$ is a symmetric integer $k\times k$ matrix,  and 
$$x_i= x^{n_i} a^{a_i} q^{q_i},$$
with $\mathbf{a}=(a_1,\ldots,a_k)$ and $\mathbf{q}=(q_1,\ldots,q_k)$ being $k$-tuples of integers, and $\mathbf{n}=(n_1,\ldots,n_k)$, $\mathbf{l}=(l_1,\ldots,l_k)$ and $\mathbf{r}=(r_1,\ldots,r_k)$ are $k$-tuples of non-negative integers, such that $\max\{l_i,r_i\}\le n_i$, $i=1,\ldots,k$.

In other words, we write the reduced $S^r$-colored HOMFLY-PT polynomial of a knot, for any $r\in \N$ as:

 \begin{eqnarray}
P_r(K)(a,q) &=& \sum_{{\tiny{\begin{array}{c}d_1,\ldots,d_k\ge 0\\
\sum_i {n_i d_i=r}\\
\end{array}}}}  \frac{(q;q)_r}{\prod_{i=1}^k {(q;q)_{d_i}}}  (-1)^{\sum_i C_{ii} d_i} a^{\sum_i a_i d_i} q^{\sum_i q_i d_i} q^{\frac{1}{2}(\sum_{i,j} C_{ij}d_id_j-\sum_i C_{ii} d_i)}  \times \nonumber\\
&&\quad\quad\quad\quad\quad\quad\times(a^{-1}q;q^{-1})_{\sum_i r_i d_i} (a^{-1}q^{-r};q^{-1})_{\sum_i l_i d_i}.\label{specred}
\end{eqnarray}

In this way, we collect all the information about all $S^r$-colored HOMFLY-PT polynomials in one symmetric integral matrix $C$, and five additional finite sequences of  integers.

\subsection{Rewriting ansatz as true quiver generating series}\label{rewrite}
In order to transform the compact forms (\ref{specform}) and (\ref{specred}) to the actual  quiver forms (\ref{kq2}), we need to perform some straightforward transformations on the additional $q$-Pochhammers. To that end we shall use couple of lemmas to simplify the expressions involving $q$-Pochhammer symbols and $q$-binomials.
\begin{lemma}[Lemma 4.5 \cite{KRSSlong}]\label{lemlong}
For any $d_1,\ldots,d_k\ge 0$, we have:
\begin{equation*}\frac{(x;q)_{d_1+\ldots+d_k}}{(q;q)_{d_1}\cdots (q;q)_{d_k}}
= \!\!\!  
\sum\limits_{{\tiny{\begin{array}{c}
\alpha_1+\beta_1=d_1\\\cdots\\\alpha_k+\beta_k=d_k
\end{array}}}}
\frac{(-x q^{-\frac{1}{2}})^{{\alpha_1+\ldots+\alpha_k}} q^{\frac{1}{2}(\alpha_1^2+\ldots+\alpha_k^2)+\sum_{i=1}^{k-1} \alpha_{i+1} (d_1+\ldots+d_i)}}{(q;q)_{\alpha_1}\cdots(q;q)_{\alpha_k}(q;q)_{\beta_1}\cdots(q;q)_{\beta_k}}      
\end{equation*}
\end{lemma}

This follows from the \textit{$q$-binomial identity:}
\begin{lemma}\label{mainlemma} We have
$$\qp{x}{k}=\sum_{i=0}^k (-1)^i x^{i} q^{\frac{1}{2}(i^2 -i)} {k \brack i} = \sum_{i=0}^k (-1)^i x^{i} q^{\frac{1}{2}(i^2 -i)} \frac{\qpp{k}}{\qpp{i}\qpp{k-i}}.$$
\end{lemma}

Another useful tool for manipulating generating series is the following lemma. In the notation below we use the quantum binomial and multinomial coefficients given by:
$${n \brack k}:=\frac{(q;q)_n}{(q;q)_k (q;q)_{n-k}}\quad\quad\textrm{ and }\quad\quad {n \brack a_1,\ldots,a_m}:=\frac{(q;q)_n}{\prod_{i=1}^m (q;q)_{a_i}},$$
where $n=\sum_{i=1}^m a_i$.
\begin{lemma}[Lemma 4.6 \cite{KRSSlong}]\label{lemvlong}
Let $a_1,\ldots,a_m$ and $b_1,\ldots,b_p$ be nonnegative integers satisfying:
$$a_1+\ldots+a_m=b_1+\ldots+b_p.$$
Then we have:
\begin{equation}\label{forexp}
\!\!\!{a_1+\ldots+a_m \brack b_1,\ldots,b_p}=
\sum_{
\begin{array}{c}
\{j_{\alpha,\beta}\}\\ 
{\scriptstyle{{\alpha=1,\ldots,p}}} \\
{\scriptstyle{\beta=1,\ldots,m}}
\end{array}
}
q^{X(\underline{j})}{a_1 \brack j_{1,1},j_{2,1},\ldots,j_{p,1}}
\cdots {a_m \brack j_{1,m},j_{2,m},\ldots,j_{p,m}},
\end{equation}

\noindent with 
\begin{equation}
\label{eqn:triangcorr} X(\underline{j})=\sum_{\scriptstyle{1\le l_1<l_2\le p}}\,\,\sum_{\scriptstyle{1\le u_1<u_2\le m}} \,j_{l_1,u_1} j_{l_2,u_2},
\end{equation}
where in formula (\ref{forexp}) we are summing over $m p$ summation indices (nonnegative integers) $j_{\alpha,\beta}$, with $\alpha=1,\ldots,p$ and $\beta=1,\ldots,m$, such that
\begin{eqnarray}
&j_{1,v}+ j_{2,v}+\ldots+ j_{p,v}=a_v,\quad v=1,\ldots,m,\\
&j_{w,1}+j_{w,2}+\ldots+j_{w,m}=b_w,\quad w=1,\ldots,p.
\end{eqnarray}
\end{lemma}

\vskip 0.3cm

Now, note that from (\ref{specred}) one can easily get the actual quiver form (\ref{kq2}), simply by expanding the two extra $q$-Pochhammers from the expression (3), by using quantum binomial formula:
\begin{equation*}
(y;q^{-1})_l=(y q^{1-l};q)_l=\sum_{i+j=l}(-1)^i y^i q^{(1-l)i}q^{\frac{1}{2}(i^2-i)} \frac{(q;q)_l}{(q;q)_i (q;q)_j}.
\end{equation*}
Moreover we can expand the $q$-Pochhammers in the following way:
\begin{equation}\label{order}(y;q)_{\sum_i l_i d_i}=(y;q)_{l_1 d_1} (y q^{l_1 d_1};q)_{l_2 d_2}\ldots (y q^{l_1 d_1+l_2 d_2+\ldots+l_{k-1} d_{k-1}};q)_{l_k d_k}.
\end{equation}
Also, for every $z$ and every $l,d\ge 1$,
\begin{eqnarray*}
(z;q)_{l d}&=&(z;q)_{d} (z q^{d};q)_{d} \ldots (z q^{(l-1) d};q)_{d}=\\
&=&\sum_{i_1+j_1=d} \sum_{i_2+j_2=d}\cdots \sum_{i_{l}+j_{l}=d}  (-1)^{\sum i_{\alpha}} z^{\sharp} q^{\sharp\sharp}{d \brack i_1}{d \brack i_2}\cdots{d \brack i_l},
\end{eqnarray*}
where by $\sharp$ and $\sharp\sharp$ we have denoted the expressions that are linear and quadratic, respectively, in $i_{\alpha}$'s and $j_{\alpha}$'s, and can be obtained by quantum binomial identity. Now, since $i_{\alpha}+j_{\alpha}=d$, for every $\alpha=1,\ldots,l$, we can use Lemma \ref{lemvlong} below, $l-1$ times, and obtain that 
$${d \brack i_1}{d \brack i_2}\cdots{d \brack i_l}={d \brack i_1,j_1}{d \brack i_2,j_2}\cdots{d \brack i_l,j_l}=\sum_{v_1,\ldots,v_{2^l}}q^{\sharp\sharp}{d\brack v_1,v_2,\ldots,v_{2^l}},$$
where the sum is over $2^l$ variables $v_i$ that sum up to $d$. 

Doing the same thing for both extra $q$-Pochhammers in (\ref{specform}), we can re-write their product as
$$(a^{-1}q;q^{-1})_{\sum_i r_i d_i} (a^{-1}q^{-\sum_i n_i d_i};q^{-1})_{\sum_i l_i d_i}=\prod_{i=1}^k \sum_{v^{(i)}_1,\ldots,v^{(i)}_{2^{l_i+r_i}}}\frac{(q;q)_{d_i}}{\prod_j (q;q)_{v^{(i)}_{j}}}.$$
Multiplying this with $\prod_{i=1}^k\frac{1}{(q;q)_{d_i}}$, gives quiver generating series in a desired form. 
In particular, the total size of the quiver is 
$$\sum_{i=1}^k 2^{l_i+r_i}.$$

\subsection{Non-uniqueness}\label{nonun}

The non-uniqueness of a quiver corresponding to a given knot, has been well known from the first results on the knots-quivers correspondence. It can come in different flavours, first that two quivers can have different sizes, but also there can be many different quivers of the same size (even minimal size) that produce the same quiver generating series, see e.g. \cite{permutahedron, KRSSlong, SW1}. In particular in \cite{permutahedron}, a detailed study of the quivers of the minimal size that correspond to a same knot was done. A discrete group of symmetries that can be represented by different and complicated permutahedra was obtained. Also, it was noticed that for the knots studied in \cite{permutahedron} the equivalence between the quivers could be obtained from the structure of differentials on HOMFLY-PT homology of the corresponding knot, \cite{DGR,GGS,GS}.   

The special form of the quiver that we propose in this paper automatically addresses large part of these equivalences. In particular the two extra $q$-Pochhammers  in (\ref{specform}) are directly motivated by the degrees of the predicted colored differentials from \cite{GGS,GS}. On the other hand, to transform this special form into a pure quiver form, as explained in Section \ref{rewrite}, depends on the order, or choice, of the terms on the RHS in (\ref{order}). Any such choice will give a different quiver with the same quiver generating series, and in particular will include a large part of the equivalences of quivers from \cite{permutahedron}. We will show this on the simple example of the knot $\mathbf{5_1}$ in Appendix \ref{ap}. 

We note that in this way we obtain different, equivalent quivers of the same size, and all of their quiver matrices have the same diagonal.

\section{Explicit results}\label{results}
In this section we present explicit results for various knots and classes of knots. We present results in the special quiver form (\ref{specform}), that is for each knot $K$ below we give a symmetric integral matrix $C$, and integral sequences ${\bf n}$, ${\bf a}$, ${\bf q}$, ${\bf l}$ and ${\bf r}$, so that the reduced $S^r$ colored HOMFLY-PT polynomial of the knot $K$ can be written as:
 \begin{eqnarray}
P_r(K)(a,q) &=& \sum_{{\tiny{\begin{array}{c}d_1,\ldots,d_k\ge 0\\
\sum_i {n_i d_i=r}\\
\end{array}}}}  \frac{(q;q)_r}{\prod_{i=1}^k {(q;q)_{d_i}}}  (-1)^{\sum_i C_{ii} d_i} a^{\sum_i a_i d_i} q^{\sum_i q_i d_i} q^{\frac{1}{2}(\sum_{i,j} C_{ij}d_id_j-\sum_i C_{ii} d_i)}  \times \nonumber\\
&&\quad\quad\quad\quad\times(a^{-1}q;q^{-1})_{\sum_i r_i d_i} (a^{-1}q^{-r};q^{-1})_{\sum_i l_i d_i}, \quad \forall r\in\N.\label{specred2}
\end{eqnarray}

\subsection{$\mathbf{9_{42}}$}\label{f942}
For knot $\bf{9_{42}}$ we have:
\begin{equation}\label{mat942}
C=\left(
\begin{array}{ccc|cc}
2&2&0&2&2\\
2&4&0&4&4\\
0&0&0&0&0\\
\hline
2&4&0&4&4\\
2&4&0&4&5
\end{array}
\right)
\end{equation}
and
$$
\mathbf{n}=\left(
\begin{array}{c}
1\\
1\\
1\\
\hline
2\\
2
\end{array}
\right), 
\quad 
\mathbf{a}=\left(
\begin{array}{c}
1\\
1\\
0\\
\hline
1\\
1\\
\end{array}
\right), 
\quad 
\mathbf{q}=\left(
\begin{array}{c}
-1\\
1\\
0\\
\hline
0\\
1
\end{array}
\right), 
\quad 
\mathbf{l}=\left(
\begin{array}{c}
1\\
1\\
0\\
\hline
2\\
2
\end{array}
\right), 
\quad 
\mathbf{r}=\left(
\begin{array}{c}
1\\
1\\
0\\
\hline
2\\
1
\end{array}
\right). 
$$
The explicit form of $P_r({\bf 9_{42}})(a,q)$ is given by:
\begin{eqnarray}
P_r({\bf 9_{42}})(a,q)&=&\sum_{d_1+d_2+d_3+2d_4+2d_5=r} \frac{(q;q)_r}{\prod_i (q;q)_{d_i}} (-1)^{d_5} a^{d_1+d_2+d_4+d_5} \times \nonumber\\
&&\quad\times\quad q^{-d_1+d_2+d_5} q^{-\frac{1}{2}(2d_1+4d_2+4d_4+5d_5)} q^{\frac{1}{2}\sum_{i,j=1}^5 C_{ij}d_id_j} \times\nonumber\\
&&\quad\times\quad
(a^{-1}q;q^{-1})_{d_1+d_2+2d_4+d_5} (a^{-1}q^{-r};q^{-1})_{d_1+d_2+2d_4+2d_5}.\label{fr942}
\end{eqnarray}

\subsection{${\bf 9_{46}}$}
For knot ${\bf 9_{46}}$ we have:
$$C=\left(
\begin{array}{ccc|c}
4&4&0&6\\
4&6&0&8\\
0&0&0&0\\
\hline
6&8&0&12
\end{array}
\right)$$
and
$$
\mathbf{n}=\left(
\begin{array}{c}
1\\
1\\
1\\
\hline
2
\end{array}
\right), 
\quad 
\mathbf{a}=\left(
\begin{array}{c}
2\\
3\\
0\\
\hline
4\\
\end{array}
\right), 
\quad 
\mathbf{q}=\left(
\begin{array}{c}
0\\
0\\
0\\
\hline
2\\
\end{array}
\right), 
\quad 
\mathbf{l}=\left(
\begin{array}{c}
1\\
1\\
0\\
\hline
2
\end{array}
\right), 
\quad 
\mathbf{r}=\left(
\begin{array}{c}
1\\
1\\
0\\
\hline
2
\end{array}
\right). 
$$
The explicit form of $P_r({\bf 9_{46}})(a,q)$ is given by:
\begin{eqnarray}
P_r({\bf 9_{46}})(a,q)&=&\sum_{d_1+d_2+d_3+2d_4=r} \frac{(q;q)_r}{\prod_i (q;q)_{d_i}}  a^{2d_1+3d_2+4d_4} \times \\
&&\quad\times\quad q^{2d_4} q^{-\frac{1}{2}(4d_1+6d_2+12d_4)} q^{\frac{1}{2}\sum_{i,j=1}^4 C_{ij}d_id_j} \times\\
&&\quad\times\quad
(a^{-1}q;q^{-1})_{d_1+d_2+2d_4} (a^{-1}q^{-r};q^{-1})_{d_1+d_2+2d_4}.
\end{eqnarray}

\subsection{$\bf 8_{20}$}
For knot ${\bf 8_{20}}$ we have:
$$C=\left(
\begin{array}{ccc|c}
0&0&0&0\\
0&3&3&4\\
0&3&5&6\\
\hline
0&4&6&8
\end{array}
\right)$$
and
$$
\mathbf{n}=\left(
\begin{array}{c}
1\\
1\\
1\\
\hline
2
\end{array}
\right), 
\quad 
\mathbf{a}=\left(
\begin{array}{c}
0\\
2\\
2\\
\hline
3\\
\end{array}
\right), 
\quad 
\mathbf{q}=\left(
\begin{array}{c}
0\\
-1\\
1\\
\hline
1\\
\end{array}
\right), 
\quad 
\mathbf{l}=\left(
\begin{array}{c}
0\\
1\\
1\\
\hline
2
\end{array}
\right), 
\quad 
\mathbf{r}=\left(
\begin{array}{c}
0\\
1\\
1\\
\hline
1
\end{array}
\right). 
$$
The explicit form of $P_r({\bf 8_{20}})(a,q)$ is given by:
\begin{eqnarray}
P_r({\bf 8_{20}})(a,q)&=&\sum_{d_1+d_2+d_3+2d_4=r} \frac{(q;q)_r}{\prod_i (q;q)_{d_i}}  (-1)^{d_2+d_3}a^{2d_2+2d_3+3d_4} \times \\
&&\quad\times\quad q^{-d_2+d_3+d_4} q^{-\frac{1}{2}(3d_2+5d_3+8d_4)} q^{\frac{1}{2}\sum_{i,j=1}^4 C_{ij}d_id_j} \times\\
&&\quad\times\quad
(a^{-1}q;q^{-1})_{d_2+d_3+d_4} (a^{-1}q^{-r};q^{-1})_{d_2+d_3+2d_4}.
\end{eqnarray}

\subsection{$\mathbf{10_{132}}$}
For knot ${\bf 10_{132}}$ we have:
$$
C=\left(\begin{array}{cccc|ccc}
5 & 5 & 1 & 3 & 6 & 8 & 8\\
5 & 7 & 1 & 3 & 8 & 10 & 10\\
1 & 1 & 0 & 1 & 1 & 2 & 2\\
3 & 3 & 1 & 3 & 5 & 6 & 6\\
\hline
6 & 8 & 1 & 5 & 10 & 12 & 13\\
8 & 10 & 2 & 6 & 12 & 16 & 16\\
8 & 10 & 2 & 6 & 13 & 16 & 17
\end{array}
\right)
$$
and
$$
\mathbf{n}=\left(
\begin{array}{c}
1\\
1\\
1\\
1\\
\hline
2\\
2\\
2
\end{array}
\right), 
\quad 
\mathbf{a}=\left(
\begin{array}{c}
3\\
3\\
1\\
2\\
\hline
4\\
5\\
5
\end{array}
\right), 
\quad 
\mathbf{q}=\left(
\begin{array}{c}
-1\\
1\\
-1\\
0\\
\hline
1\\
3\\
4
\end{array}
\right), 
\quad 
\mathbf{l}=\left(
\begin{array}{c}
1\\
1\\
0\\
0\\
\hline
2\\
2\\
2
\end{array}
\right), 
\quad 
\mathbf{r}=\left(
\begin{array}{c}
1\\
1\\
0\\
1\\
\hline
1\\
2\\
2
\end{array}
\right). 
$$
The explicit form of $P_r({\bf 10_{132}})(a,q)$ is given by:
\begin{eqnarray*}
P_r({\bf 10_{132}})(a,q)&=&\sum_{d_1+d_2+d_3+d_4+2d_5+2d_6+2d_7=r} \frac{(q;q)_r}{\prod_i (q;q)_{d_i}} \,\,(-1)^{d_1+d_2+d_4+d_7} \times \\
&&\quad\times\quad a^{3d_1+3d_2+d_3+2d_4+4d_5+5d_6+5d_7} q^{-d_1+d_2-d_3+d_5+3d_6+4d_7} \times \\
&&\quad \times\quad q^{-\frac{1}{2}(5d_1+7d_2+3d_4+10d_5+16d_6+17d_7)} q^{\frac{1}{2}\sum_{i,j=1}^7 C_{ij}d_id_j}\times\\
&&\quad\times \quad
(a^{-1}q;q^{-1})_{d_1+d_2+d_4+d_5+2d_6+2d_7} (a^{-1}q^{-r};q^{-1})_{d_1+d_2+2d_5+2d_6+2d_7}.
\end{eqnarray*}

\subsection{${\bf 10_{145}}$}
For knot ${\bf 10_{145}}$ we have:
$$
C=\left(\begin{array}{ccccc|ccc}
9 & 7 & 2 & 4 & 5 & 13 & 14 & 14\\
7 & 7 & 2 & 4 & 5 & 11 & 12 & 12\\
2 & 2 & 0 & 1 & 3 & 3 & 4 & 5\\
4 & 4 & 1 & 3 & 4 & 7 & 8 & 8\\
5 & 5 & 3 & 4 & 5 & 10 & 10 & 10\\
\hline
13 & 11 & 3 & 7 & 10 & 21 & 22 & 23\\
14 & 12 & 4 & 8 & 10 & 22 & 24 & 24\\
14 & 12 & 5 & 8 & 10 & 23 & 24 & 25
\end{array}
\right)
$$
and
$$
\mathbf{n}=\left(
\begin{array}{c}
1\\
1\\
1\\
1\\
1\\
\hline
2\\
2\\
2
\end{array}
\right), 
\quad 
\mathbf{a}=\left(
\begin{array}{c}
5\\
4\\
2\\
3\\
3\\
\hline
8\\
8\\
8
\end{array}
\right), 
\quad 
\mathbf{q}=\left(
\begin{array}{c}
0\\
0\\
-2\\
-1\\
1\\
\hline
3\\
5\\
6
\end{array}
\right), 
\quad 
\mathbf{l}=\left(
\begin{array}{c}
1\\
1\\
0\\
0\\
0\\
\hline
2\\
2\\
2
\end{array}
\right), 
\quad 
\mathbf{r}=\left(
\begin{array}{c}
1\\
1\\
0\\
1\\
1\\
\hline
2\\
2\\
2
\end{array}
\right). 
$$
The explicit form of $P_r({\bf 10_{145}})(a,q)$ is given by:
\begin{eqnarray*}
P_r({\bf 10_{145}})(a,q)&=&\sum_{d_1+d_2+d_3+d_4+d_5+2d_6+2d_7+2d_8=r} \frac{(q;q)_r}{\prod_i (q;q)_{d_i}} \,\,(-1)^{d_1+d_2+d_4+d_5+d_6+d_8} \times \\
&&\quad\times\quad a^{5d_1+4d_2+2d_3+3d_4+3d_5+8d_6+8d_7+8d_8} q^{-2d_3-d_4+d_5+3d_6+5d_7+6d_8} \times \\
&&\quad \times\quad q^{-\frac{1}{2}(9d_1+7d_2+3d_4+5d_5+21d_6+24d_7+25d_8)} q^{\frac{1}{2}\sum_{i,j=1}^8 C_{ij}d_id_j}\times\\
&&\quad\times \quad
(a^{-1}q;q^{-1})_{d_1+d_2+d_4+d_5+2d_6+2d_7+2d_8} (a^{-1}q^{-r};q^{-1})_{d_1+d_2+2d_6+2d_7+2d_8}.
\end{eqnarray*}

\subsection{${\bf 10_{139}}$}
For knot ${\bf 10_{139}}$ we have:
$$
C=\left(\begin{array}{cccccccc}
0 & 1 & 3 & 5 & 7 & 3 & 5 & 4\\
1 & 3 & 4 & 6 & 8 & 5 & 6 & 6\\
3 & 4 & 5 & 6 & 8 & 7 & 7 & 7\\
5 & 6 & 6 & 7 & 8 & 8 & 8 & 8\\
7 & 8 & 8 & 8 & 9 & 9 & 9 & 9\\
3 & 5 & 7 & 8 & 9 & 8 & 9 & 8\\
5 & 6 & 7 & 8 & 9 & 9 & 10 & 9\\
4 & 6 & 7 & 8 & 9 & 8 & 9 & 9
\end{array}
\right)
$$
and
$$
\mathbf{n}=\left(
\begin{array}{c}
1\\
1\\
1\\
1\\
1\\
1\\
1\\
1
\end{array}
\right), 
\quad 
\mathbf{a}=\left(
\begin{array}{c}
4\\
5\\
5\\
5\\
5\\
6\\
6\\
6
\end{array}
\right), 
\quad 
\mathbf{q}=\left(
\begin{array}{c}
-4\\
-3\\
-1\\
1\\
3\\
-1\\
1\\
0
\end{array}
\right), 
\quad 
\mathbf{l}=\left(
\begin{array}{c}
0\\
0\\
0\\
0\\
0\\
1\\
1\\
1
\end{array}
\right), 
\quad 
\mathbf{r}=\left(
\begin{array}{c}
0\\
1\\
1\\
1\\
1\\
1\\
1\\
1
\end{array}
\right). 
$$

\subsection{${\bf 10_{152}}$}
For knot ${\bf 10_{152}}$ we have:
$$
C=\left(\begin{array}{cccccccccc|c}
0 & 1 & 3 & 5 & 7 & 3 & 5 & 3 & 4 & 5 & 9\\
1 & 3 & 4 & 6 & 8 & 5 & 6 & 5 & 6 & 6 & 14\\
3 & 4 & 5 & 6 & 8 & 7 & 7 & 7 & 7 & 7 & 18\\
5 & 6 & 6 & 7 & 8 & 8 & 8 & 8 & 8 & 8 & 22\\
7 & 8 & 8 & 8 & 9 & 9 & 9 & 9 & 9 & 9 & 26\\
3 & 5 & 7 & 8 & 9 & 8 & 9 & 7 & 8 & 8 & 21\\
5 & 6 & 7 & 8 & 9 & 9 & 10 & 8 &9 & 9 & 24\\
3 & 5 & 7 & 8 & 9 & 7 & 8 & 8 & 8 & 9 & 21\\
4 & 6 & 7 & 8 & 9 & 8 & 9 & 8 & 9 & 9 & 23\\
5 & 6 & 7 & 8 & 9 & 8 & 9 & 9 &9 & 10 & 24\\
\hline
9 & 14 & 18 & 22 & 26 & 21 & 24 & 21 &23 & 24 & 60
\end{array}
\right)
$$
and
$$
\mathbf{n}=\left(
\begin{array}{c}
1\\
1\\
1\\
1\\
1\\
1\\
1\\
1\\
1\\
1\\
\hline
3
\end{array}
\right), 
\quad 
\mathbf{a}=\left(
\begin{array}{c}
4\\
5\\
5\\
5\\
5\\
6\\
6\\
6\\
6\\
6\\
\hline
17
\end{array}
\right), 
\quad 
\mathbf{q}=\left(
\begin{array}{c}
-4\\
-3\\
-1\\
1\\
3\\
-1\\
1\\
-1\\
0\\
1\\
\hline
16
\end{array}
\right), 
\quad 
\mathbf{l}=\left(
\begin{array}{c}
0\\
0\\
0\\
0\\
0\\
1\\
1\\
1\\
1\\
1\\
\hline
2
\end{array}
\right), 
\quad 
\mathbf{r}=\left(
\begin{array}{c}
0\\
1\\
1\\
1\\
1\\
1\\
1\\
1\\
1\\
1\\
\hline
3
\end{array}
\right). 
$$

\subsection{$(2,2p+1)$ torus knots}
For $T_{2,2p+1}$ torus knots, for arbitrary $p\ge 1$, the matrix $C$ is the following  $(p+1)\times (p+1)$:
 
$$C=\left(
\begin{array}{ccccccc}
0&1&3&5&\cdots&2p-3&2p-1\\
1&3&4&6&\cdots&2p-2&2p\\
3&4&5&6&\cdots&2p-2&2p\\
5&6&6&7&\cdots&2p-2&2p\\
\vdots&\vdots&\vdots&\vdots&\ddots&\vdots&\vdots\\
2p-3&2p-2&2p-2&2p-2&\cdots&2p-1&2p\\
2p-1&2p&2p&2p&\cdots&2p&2p+1
\end{array}
\right),$$
while the corresponding vectors are of length $p+1$, and are given by
$$
\mathbf{n}=\left(
\begin{array}{c}
1\\
1\\
1\\
1\\
\vdots\\
1
\end{array}
\right), 
\quad 
\mathbf{a}=\left(
\begin{array}{c}
p\\
p+1\\
p+1\\
p+1\\
\vdots\\
p+1
\end{array}
\right), 
\quad 
\mathbf{q}=\left(
\begin{array}{c}
-p\\
1-p\\
3-p\\
5-p\\
\vdots\\
p-1
\end{array}
\right), 
\quad 
\mathbf{l}=\left(
\begin{array}{c}
0\\
0\\
0\\
0\\
\vdots\\
0
\end{array}
\right), 
\quad 
\mathbf{r}=\left(
\begin{array}{c}
0\\
1\\
1\\
1\\
\vdots\\
1
\end{array}
\right). 
$$ 
 
Explicitely, for couple of first $(2,2p+1)$ torus knots, we have the following:\\

\noindent{$\bullet \quad \mathbf{3_1}$} \indent For knot ${\bf 3_{1}}$ we have:
$$C=\left(
\begin{array}{cc}
0&1\\
1&3
\end{array}
\right)$$
and
$$
\mathbf{n}=\left(
\begin{array}{c}
1\\
1
\end{array}
\right), 
\quad 
\mathbf{a}=\left(
\begin{array}{c}
1\\
2
\end{array}
\right), 
\quad 
\mathbf{q}=\left(
\begin{array}{c}
-1\\
0
\end{array}
\right), 
\quad 
\mathbf{l}=\left(
\begin{array}{c}
0\\
0
\end{array}
\right), 
\quad 
\mathbf{r}=\left(
\begin{array}{c}
0\\
1\end{array}
\right). 
$$

\noindent
{$\bullet \quad \mathbf{5_1}$} \indent For knot ${\bf 5_{1}}$ we have:
\begin{equation}\label{mat51}
C=\left(
\begin{array}{ccc}
0&1&3\\
1&3&4\\
3&4&5
\end{array}
\right)
\end{equation}
and
$$
\mathbf{n}=\left(
\begin{array}{c}
1\\
1\\
1
\end{array}
\right), 
\quad 
\mathbf{a}=\left(
\begin{array}{c}
2\\
3\\
3
\end{array}
\right), 
\quad 
\mathbf{q}=\left(
\begin{array}{c}
-2\\
-1\\
1
\end{array}
\right), 
\quad 
\mathbf{l}=\left(
\begin{array}{c}
0\\
0\\
0
\end{array}
\right), 
\quad 
\mathbf{r}=\left(
\begin{array}{c}
0\\
1\\
1
\end{array}
\right). 
$$

\noindent
{$\bullet \quad \mathbf{7_1}$} \indent For knot ${\bf 7_{1}}$ we have:
$$C=\left(
\begin{array}{cccc}
0&1&3&5\\
1&3&4&6\\
3&4&5&6\\
5&6&6&7
\end{array}
\right)$$
and
$$
\mathbf{n}=\left(
\begin{array}{c}
1\\
1\\
1\\
1
\end{array}
\right), 
\quad 
\mathbf{a}=\left(
\begin{array}{c}
3\\
4\\
4\\
4
\end{array}
\right), 
\quad 
\mathbf{q}=\left(
\begin{array}{c}
-3\\
-2\\
0\\
2
\end{array}
\right), 
\quad 
\mathbf{l}=\left(
\begin{array}{c}
0\\
0\\
0\\
0
\end{array}
\right), 
\quad 
\mathbf{r}=\left(
\begin{array}{c}
0\\
1\\
1\\
1
\end{array}
\right). 
$$

\noindent
{$\bullet \quad \mathbf{9_1}$} \indent For knot ${\bf 9_{1}}$ we have:
$$C=\left(
\begin{array}{ccccc}
0&1&3&5&7\\
1&3&4&6&8\\
3&4&5&6&8\\
5&6&6&7&8\\
7&8&8&8&9
\end{array}
\right)$$
and
$$
\mathbf{n}=\left(
\begin{array}{c}
1\\
1\\
1\\
1\\
1
\end{array}
\right), 
\quad 
\mathbf{a}=\left(
\begin{array}{c}
4\\
5\\
5\\
5\\
5
\end{array}
\right), 
\quad 
\mathbf{q}=\left(
\begin{array}{c}
-4\\
-3\\
-1\\
1\\
3
\end{array}
\right), 
\quad 
\mathbf{l}=\left(
\begin{array}{c}
0\\
0\\
0\\
0\\
0
\end{array}
\right), 
\quad 
\mathbf{r}=\left(
\begin{array}{c}
0\\
1\\
1\\
1\\
1
\end{array}
\right). 
$$


\subsection{Twist knots I -- even series}
In this section we give results for the twist knots with even numbers of crossings.
Let $TK_{2p+2}$, with $p\ge 1$ be a twist knots with $2p+2$ crossings. Then the corresponding data is given by the following $(p+1)\times (p+1)$ matrix:

$$C=\left(
\begin{array}{cccccc}
0&0&0&\cdots&0&0\\
0&2&2&\cdots&2&2\\
0&2&4&\cdots&4&4\\
\vdots&\vdots&\vdots&\ddots&\vdots&\vdots\\
0&2&4&\cdots&2p-2&2p-2\\
0&2&4&\cdots&2p-2&2p
\end{array}
\right)$$
and the vectors of length $p+1$
$$
\mathbf{n}=\left(
\begin{array}{c}
1\\
1\\
1\\
\vdots\\
1
\end{array}
\right), 
\quad 
\mathbf{a}=\left(
\begin{array}{c}
0\\
1\\
2\\
\vdots\\
p
\end{array}
\right), 
\quad 
\mathbf{q}=\left(
\begin{array}{c}
0\\
0\\
0\\
\vdots\\
0
\end{array}
\right), 
\quad 
\mathbf{l}=\left(
\begin{array}{c}
0\\
1\\
1\\
\vdots\\
1
\end{array}
\right), 
\quad 
\mathbf{r}=\left(
\begin{array}{c}
0\\
1\\
1\\
\vdots\\
1
\end{array}
\right). 
$$
 
Explicitely, for the first couple of twist knots from even series we have:\\

\noindent {$\bullet \quad \mathbf{4_1}$} \indent For knot ${\bf 4_{1}}$ we have:
$$C=\left(
\begin{array}{cc}
0&0\\
0&2
\end{array}
\right)$$
and
$$
\mathbf{n}=\left(
\begin{array}{c}
1\\
1
\end{array}
\right), 
\quad 
\mathbf{a}=\left(
\begin{array}{c}
0\\
1
\end{array}
\right), 
\quad 
\mathbf{q}=\left(
\begin{array}{c}
0\\
0
\end{array}
\right), 
\quad 
\mathbf{l}=\left(
\begin{array}{c}
0\\
1
\end{array}
\right), 
\quad 
\mathbf{r}=\left(
\begin{array}{c}
0\\
1\end{array}
\right). 
$$

\noindent
{$\bullet \quad \mathbf{6_1}$} \indent For knot ${\bf 6_{1}}$ we have:
$$C=\left(
\begin{array}{ccc}
0&0&0\\
0&2&2\\
0&2&4
\end{array}
\right)$$
and
$$
\mathbf{n}=\left(
\begin{array}{c}
1\\
1\\
1
\end{array}
\right), 
\quad 
\mathbf{a}=\left(
\begin{array}{c}
0\\
1\\
2
\end{array}
\right), 
\quad 
\mathbf{q}=\left(
\begin{array}{c}
0\\
0\\
0
\end{array}
\right), 
\quad 
\mathbf{l}=\left(
\begin{array}{c}
0\\
1\\
1
\end{array}
\right), 
\quad 
\mathbf{r}=\left(
\begin{array}{c}
0\\
1\\
1
\end{array}
\right). 
$$

\noindent
{$\bullet \quad \mathbf{8_1}$} \indent For knot ${\bf 8_{1}}$ we have:
$$C=\left(
\begin{array}{cccc}
0&0&0&0\\
0&2&2&2\\
0&2&4&4\\
0&2&4&6
\end{array}
\right)$$
and
$$
\mathbf{n}=\left(
\begin{array}{c}
1\\
1\\
1\\
1
\end{array}
\right), 
\quad 
\mathbf{a}=\left(
\begin{array}{c}
0\\
1\\
2\\
3
\end{array}
\right), 
\quad 
\mathbf{q}=\left(
\begin{array}{c}
0\\
0\\
0\\
0
\end{array}
\right), 
\quad 
\mathbf{l}=\left(
\begin{array}{c}
0\\
1\\
1\\
1
\end{array}
\right), 
\quad 
\mathbf{r}=\left(
\begin{array}{c}
0\\
1\\
1\\
1
\end{array}
\right). 
$$

\noindent
{$\bullet \quad \mathbf{10_1}$} \indent For knot ${\bf 10_{1}}$ we have:
$$C=\left(
\begin{array}{ccccc}
0&0&0&0&0\\
0&2&2&2&2\\
0&2&4&4&4\\
0&2&4&6&6\\
0&2&4&6&8
\end{array}
\right)$$
and
$$
\mathbf{n}=\left(
\begin{array}{c}
1\\
1\\
1\\
1\\
1
\end{array}
\right), 
\quad 
\mathbf{a}=\left(
\begin{array}{c}
0\\
1\\
2\\
3\\
4
\end{array}
\right), 
\quad 
\mathbf{q}=\left(
\begin{array}{c}
0\\
0\\
0\\
0\\
0
\end{array}
\right), 
\quad 
\mathbf{l}=\left(
\begin{array}{c}
0\\
1\\
1\\
1\\
1
\end{array}
\right), 
\quad 
\mathbf{r}=\left(
\begin{array}{c}
0\\
1\\
1\\
1\\
1
\end{array}
\right). 
$$


\subsection{Twist knots II -- odd series}
In this section we give results for the twist knots with odd numbers of crossings.
Let $TK_{2p+1}$, with $p\ge 1$ be a twist knots with $2p+1$ crossings. Then the corresponding data is given by the following $(p+1)\times (p+1)$ matrix:

$$C=\left(
\begin{array}{cccccc}
0&1&1&\cdots&1&1\\
1&3&3&\cdots&3&3\\
1&3&5&\cdots&5&5\\
\vdots&\vdots&\vdots&\ddots&\vdots&\vdots\\
1&3&5&\cdots&2p-1&2p-1\\
1&3&5&\cdots&2p-1&2p+1
\end{array}
\right)$$
and the vectors of length $p+1$
$$
\mathbf{n}=\left(
\begin{array}{c}
1\\
1\\
1\\
\vdots\\
1
\end{array}
\right), 
\quad 
\mathbf{a}=\left(
\begin{array}{c}
1\\
2\\
3\\
\vdots\\
p+1
\end{array}
\right), 
\quad 
\mathbf{q}=\left(
\begin{array}{c}
-1\\
0\\
0\\
\vdots\\
0
\end{array}
\right), 
\quad 
\mathbf{l}=\left(
\begin{array}{c}
0\\
0\\
1\\
\vdots\\
1
\end{array}
\right), 
\quad 
\mathbf{r}=\left(
\begin{array}{c}
0\\
1\\
1\\
\vdots\\
1
\end{array}
\right). 
$$

The first knot in the series is ${\bf 3_1}$, for which we have already written the explicit form in the 2-strand torus knots series, and it agrees with the above general formula for $p=1$. As for couple of next twist knots with odd number of crossings, the explicit expressions are the following:\\

\noindent{$\bullet \quad \mathbf{5_2}$} \indent For knot ${\bf 5_{2}}$ we have:
$$C=\left(
\begin{array}{ccc}
0&1&1\\
1&3&3\\
1&3&5
\end{array}
\right)$$
and
$$
\mathbf{n}=\left(
\begin{array}{c}
1\\
1\\
1
\end{array}
\right), 
\quad 
\mathbf{a}=\left(
\begin{array}{c}
1\\
2\\
3
\end{array}
\right), 
\quad 
\mathbf{q}=\left(
\begin{array}{c}
-1\\
0\\
0
\end{array}
\right), 
\quad 
\mathbf{l}=\left(
\begin{array}{c}
0\\
0\\
1
\end{array}
\right), 
\quad 
\mathbf{r}=\left(
\begin{array}{c}
0\\
1\\
1
\end{array}
\right). 
$$

\noindent
{$\bullet \quad \mathbf{7_2}$} \indent For knot ${\bf 7_{2}}$ we have:
$$C=\left(
\begin{array}{cccc}
0&1&1&1\\
1&3&3&3\\
1&3&5&5\\
1&3&5&7
\end{array}
\right)$$
and
$$
\mathbf{n}=\left(
\begin{array}{c}
1\\
1\\
1\\
1
\end{array}
\right), 
\quad 
\mathbf{a}=\left(
\begin{array}{c}
1\\
2\\
3\\
4
\end{array}
\right), 
\quad 
\mathbf{q}=\left(
\begin{array}{c}
-1\\
0\\
0\\
0
\end{array}
\right), 
\quad 
\mathbf{l}=\left(
\begin{array}{c}
0\\
0\\
1\\
1
\end{array}
\right), 
\quad 
\mathbf{r}=\left(
\begin{array}{c}
0\\
1\\
1\\
1
\end{array}
\right). 
$$

\noindent
{$\bullet \quad \mathbf{9_2}$} \indent For knot ${\bf 9_{2}}$ we have:
$$C=\left(
\begin{array}{ccccc}
0&1&1&1&1\\
1&3&3&3&3\\
1&3&5&5&5\\
1&3&5&7&7\\
1&3&5&7&9
\end{array}
\right)$$
and
$$
\mathbf{n}=\left(
\begin{array}{c}
1\\
1\\
1\\
1\\
1
\end{array}
\right), 
\quad 
\mathbf{a}=\left(
\begin{array}{c}
1\\
2\\
3\\
4\\
5
\end{array}
\right), 
\quad 
\mathbf{q}=\left(
\begin{array}{c}
-1\\
0\\
0\\
0\\
0
\end{array}
\right), 
\quad 
\mathbf{l}=\left(
\begin{array}{c}
0\\
0\\
1\\
1\\
1
\end{array}
\right), 
\quad 
\mathbf{r}=\left(
\begin{array}{c}
0\\
1\\
1\\
1\\
1
\end{array}
\right). 
$$

\subsection{Some other small knots with up to 7 crossings}

\noindent
{$\bullet \quad \mathbf{6_2}$} \indent For knot ${\bf 6_{2}}$ we have:
$$C=\left(
\begin{array}{cccc}
0&1&0&2\\
1&3&2&3\\
0&2&2&3\\
2&3&3&4
\end{array}
\right)$$
and
$$
\mathbf{n}=\left(
\begin{array}{c}
1\\
1\\
1\\
1
\end{array}
\right), 
\quad 
\mathbf{a}=\left(
\begin{array}{c}
1\\
2\\
2\\
2
\end{array}
\right), 
\quad 
\mathbf{q}=\left(
\begin{array}{c}
-1\\
0\\
-1\\
1
\end{array}
\right), 
\quad 
\mathbf{l}=\left(
\begin{array}{c}
0\\
0\\
1\\
1
\end{array}
\right), 
\quad 
\mathbf{r}=\left(
\begin{array}{c}
0\\
1\\
1\\
1
\end{array}
\right). 
$$

\noindent
{$\bullet \quad \mathbf{7_3}$} \indent For knot ${\bf 7_{3}}$ we have:
$$C=\left(
\begin{array}{ccccc}
0&1&3&1&3\\
1&3&4&3&4\\
3&4&5&5&5\\
1&3&5&5&6\\
3&4&5&6&7
\end{array}
\right)$$
and
$$
\mathbf{n}=\left(
\begin{array}{c}
1\\
1\\
1\\
1\\
1
\end{array}
\right), 
\quad 
\mathbf{a}=\left(
\begin{array}{c}
2\\
3\\
3\\
4\\
4
\end{array}
\right), 
\quad 
\mathbf{q}=\left(
\begin{array}{c}
-2\\
-1\\
1\\
-1\\
1
\end{array}
\right), 
\quad 
\mathbf{l}=\left(
\begin{array}{c}
0\\
0\\
0\\
1\\
1
\end{array}
\right), 
\quad 
\mathbf{r}=\left(
\begin{array}{c}
0\\
1\\
1\\
1\\
1
\end{array}
\right). 
$$

\noindent
{$\bullet \quad \mathbf{7_4}$} \indent For knot ${\bf 7_{4}}$ we have:
$$C=\left(
\begin{array}{ccccc}
0&1&1&1&1\\
1&3&3&3&3\\
1&3&5&4&5\\
1&3&4&5&5\\
1&3&5&5&7
\end{array}
\right)$$
and
$$
\mathbf{n}=\left(
\begin{array}{c}
1\\
1\\
1\\
1\\
1
\end{array}
\right), 
\quad 
\mathbf{a}=\left(
\begin{array}{c}
1\\
2\\
3\\
3\\
4
\end{array}
\right), 
\quad 
\mathbf{q}=\left(
\begin{array}{c}
-1\\
0\\
0\\
0\\
0
\end{array}
\right), 
\quad 
\mathbf{l}=\left(
\begin{array}{c}
0\\
0\\
1\\
1\\
1
\end{array}
\right), 
\quad 
\mathbf{r}=\left(
\begin{array}{c}
0\\
1\\
1\\
1\\
1
\end{array}
\right). 
$$

\noindent
{$\bullet \quad \mathbf{7_5}$} \indent For knot ${\bf 7_{5}}$ we have:
$$C=\left(
\begin{array}{cccccc}
0&1&3&1&3&2\\
1&3&4&3&4&4\\
3&4&5&5&5&5\\
1&3&5&5&6&5\\
3&4&5&6&7&6\\
2&4&5&5&6&6
\end{array}
\right)$$
and
$$
\mathbf{n}=\left(
\begin{array}{c}
1\\
1\\
1\\
1\\
1\\
1
\end{array}
\right), 
\quad 
\mathbf{a}=\left(
\begin{array}{c}
2\\
3\\
3\\
4\\
4\\
4
\end{array}
\right), 
\quad 
\mathbf{q}=\left(
\begin{array}{c}
-2\\
-1\\
1\\
-1\\
1\\
0
\end{array}
\right), 
\quad 
\mathbf{l}=\left(
\begin{array}{c}
0\\
0\\
0\\
1\\
1\\
1
\end{array}
\right), 
\quad 
\mathbf{r}=\left(
\begin{array}{c}
0\\
1\\
1\\
1\\
1\\
1
\end{array}
\right). 
$$

\noindent
{$\bullet \quad \mathbf{7_6}$} \indent For knot ${\bf 7_{6}}$ we have:
$$C=\left(
\begin{array}{cccccc}
0&1&0&2&1&1\\
1&3&2&3&3&3\\
0&2&2&3&2&3\\
2&3&3&4&3&4\\
1&3&2&3&3&4\\
1&3&3&4&4&5
\end{array}
\right)$$
and
$$
\mathbf{n}=\left(
\begin{array}{c}
1\\
1\\
1\\
1\\
1\\
1
\end{array}
\right), 
\quad 
\mathbf{a}=\left(
\begin{array}{c}
1\\
2\\
2\\
2\\
2\\
3
\end{array}
\right), 
\quad 
\mathbf{q}=\left(
\begin{array}{c}
-1\\
0\\
-1\\
1\\
0\\
0
\end{array}
\right), 
\quad 
\mathbf{l}=\left(
\begin{array}{c}
0\\
0\\
1\\
1\\
1\\
1
\end{array}
\right), 
\quad 
\mathbf{r}=\left(
\begin{array}{c}
0\\
1\\
1\\
1\\
1\\
1
\end{array}
\right). 
$$


\subsection{Thick torus knots}
\noindent
{$\bullet \quad \mathbf{8_{19}}$ ({\bf i.e.} $\mathbf{T(3,4)})$} \indent For knot $8_{19}$ we have:
$$C=\left(
\begin{array}{ccccc}
0&1&3&5&3\\
1&3&4&6&5\\
3&4&5&6&6\\
5&6&6&7&7\\
3&5&6&7&8
\end{array}
\right)$$
and
$$
\mathbf{n}=\left(
\begin{array}{c}
1\\
1\\
1\\
1\\
1
\end{array}
\right), 
\quad 
\mathbf{a}=\left(
\begin{array}{c}
3\\
4\\
4\\
4\\
5
\end{array}
\right), 
\quad 
\mathbf{q}=\left(
\begin{array}{c}
-3\\
-2\\
0\\
2\\
0
\end{array}
\right), 
\quad 
\mathbf{l}=\left(
\begin{array}{c}
0\\
0\\
0\\
0\\
1
\end{array}
\right), 
\quad 
\mathbf{r}=\left(
\begin{array}{c}
0\\
1\\
1\\
1\\
1
\end{array}
\right). 
$$

\noindent
{$\bullet \quad \mathbf{10_{124}}$ ({\bf i.e. } $\mathbf{T(3,5)}$)} \indent For knot $10_{124}$ we have:
$$C=\left(
\begin{array}{ccccccc}
0&1&3&5&7&3&5\\
1&3&4&6&8&5&6\\
3&4&5&6&8&7&7\\
5&6&6&7&8&8&8\\
7&8&8&8&9&9&9\\
3&5&7&8&9&8&9\\
5&6&7&8&9&9&10
\end{array}
\right)$$
and
$$
\mathbf{n}=\left(
\begin{array}{c}
1\\
1\\
1\\
1\\
1\\
1\\
1
\end{array}
\right), 
\quad 
\mathbf{a}=\left(
\begin{array}{c}
4\\
5\\
5\\
5\\
5\\
6\\
6
\end{array}
\right), 
\quad 
\mathbf{q}=\left(
\begin{array}{c}
-4\\
-3\\
-1\\
1\\
3\\
-1\\
1
\end{array}
\right), 
\quad 
\mathbf{l}=\left(
\begin{array}{c}
0\\
0\\
0\\
0\\
0\\
1\\
1
\end{array}
\right), 
\quad 
\mathbf{r}=\left(
\begin{array}{c}
0\\
1\\
1\\
1\\
1\\
1\\
1
\end{array}
\right). 
$$


\subsection{Some special ones}
{$\bullet \quad \mathbf{6_{3}}$} \indent For knot ${\bf 6_3}$ there are two options: the first one has all quiver generators of level 1 (the vector $\mathbf{n}$ consists of all 1's), but it doesn't have the special form (\ref{specform}). The other option has the form (\ref{specform}) but also has quiver generators at level 2.

The higher level option of the form (\ref{specform}) is given by:

$$C=\left(
\begin{array}{cccc|c}
0&0&-1&1&0\\
0&2&1&2&3\\
-1&1&1&2&2\\
1&2&2&3&4\\
\hline
0&3&2&4&6
\end{array}
\right)$$
and
$$
\mathbf{n}=\left(
\begin{array}{c}
1\\
1\\
1\\
1\\
\hline
2
\end{array}
\right), 
\quad 
\mathbf{a}=\left(
\begin{array}{c}
0\\
1\\
1\\
1\\
\hline
2
\end{array}
\right), 
\quad 
\mathbf{q}=\left(
\begin{array}{c}
0\\
0\\
-1\\
1\\
\hline
1
\end{array}
\right), 
\quad 
\mathbf{l}=\left(
\begin{array}{c}
0\\
1\\
1\\
1\\
\hline
2
\end{array}
\right), 
\quad 
\mathbf{r}=\left(
\begin{array}{c}
0\\
1\\
1\\
1\\
\hline
1
\end{array}
\right). 
$$

The other option has slightly different form comparing to (\ref{specform}): instead of the "positive differential" part
$(a^{-1}q;q^{-1})_{\sum_i r_i d_i}$, which in this case should be given by $$(a^{-1}q;q^{-1})_{d_2+d_3+d_4}=(a^{-1}q;q^{-1})_{d_2+d_3}(a^{-1}q^{1-d_2-d_3};q^{-1})_{d_4},$$
should be replaced by slightly shifted 
$$(a^{-1}q;q^{-1})_{d_2+d_3}(a^{-1}q^{1-d_2};q^{-1})_{d_4}.$$ 
Note also that this expression is symmetric with respect to $d_3$ and $d_4$:
$$(a^{-1}q;q^{-1})_{d_2+d_3}(a^{-1}q^{1-d_2};q^{-1})_{d_4}=\frac{(a^{-1}q;q^{-1})_{d_2+d_3}(a^{-1}q;q^{-1})_{d_2+d_4}}{(a^{-1}q;q^{-1})_{d_2}}.$$
Clearly this still gives quiver form.

The quiver data is given by:
$$C=\left(
\begin{array}{cccc}
0&0&0&0\\
0&2&1&2\\
0&1&1&1\\
0&2&1&3
\end{array}
\right)$$
and
$$
\mathbf{n}=\left(
\begin{array}{c}
1\\
1\\
1\\
1
\end{array}
\right), 
\quad 
\mathbf{a}=\left(
\begin{array}{c}
0\\
1\\
1\\
1
\end{array}
\right), 
\quad 
\mathbf{q}=\left(
\begin{array}{c}
0\\
0\\
-1\\
1
\end{array}
\right), 
\quad 
\mathbf{l}=\left(
\begin{array}{c}
0\\
1\\
1\\
1
\end{array}
\right), 
\quad 
\mathbf{r}=\left(
\begin{array}{c}
0\\
1\\
1\\
1
\end{array}
\right). 
$$

\section{Knot {$\mathbf{9_{42}}$}}\label{sec942}
\subsection{Colored HOMFLY-PT polynomial}
First we shall focus on the decategorified (or unrefined) case, and look only for the $S^r$-colored HOMFLY-PT polynomials. 

The explicit formula for the reduced $S^r$-colored HOMFLY-PT polynomial of $\mathbf{9_{42}}$ is given by (\ref{fr942}). It can also be rewritten in the following form:
\begin{eqnarray}
P_r(\mathbf{9_{42}})(a,q)&=&\sum_{k+2l =r} \frac{(q;q)_r}{(q;q)_k (q;q)_l} \, \sum_{0\le i \le j \le k} \sum_{\alpha=0}^l   (-1)^{\alpha} a^{l+j} q^{i^2+j^2-2j+i+2l(i+j)} q^{2(l^2-l)} q^{\frac{1}{2}(\alpha^2+\alpha)}\times\nonumber\\
&&\quad\quad\times  \left[\!\begin{array}{c} k\\ j\end{array}\!\right]\left[\!\begin{array}{c} j\\ i\end{array}\!\right]\left[\!\begin{array}{c} l\\ \alpha\end{array}\!\right]\, (a^{-1} q; q^{-1})_{2l+j-\alpha} (a^{-1} q^{-r}; q^{-1})_{2l+j}.\label{glpol}
\end{eqnarray}

The generating series can be rewritten in the quiver form:
\begin{eqnarray}
&&P(\mathbf{9_{42}})(x;a,q)=1+ \sum_{r\ge 1}P_r (\mathbf{9_{42}})(a,q) x^r = \\
&&\quad=\sum_{\mathbf{d}=(\mathbf{d^1},\mathbf{d^2})} {(-1)^{\sum_i t_i d_i} a^{\sum_i a_i d_i} q^{\sum_i q_i d_i} q^{\frac{1}{2}\sum_{i,j}C_{i,j}d_i d_j-\frac{1}{2}\sum_i C_{i,i} d_i} \frac{(q;q)_r}{\prod_i (q;q)_{d_i}}} x^{\sum_i {d^1_i}+ 2\sum_i {d^2_i}} .\nonumber
\end{eqnarray}
In this case, the corresponding quiver has size 33. There are 9 generators of level 1, and 24 generators of level 2, and the corresponding quiver matrix is:

\begin{figure}[h]
  \centering
  \includegraphics[scale=0.55]{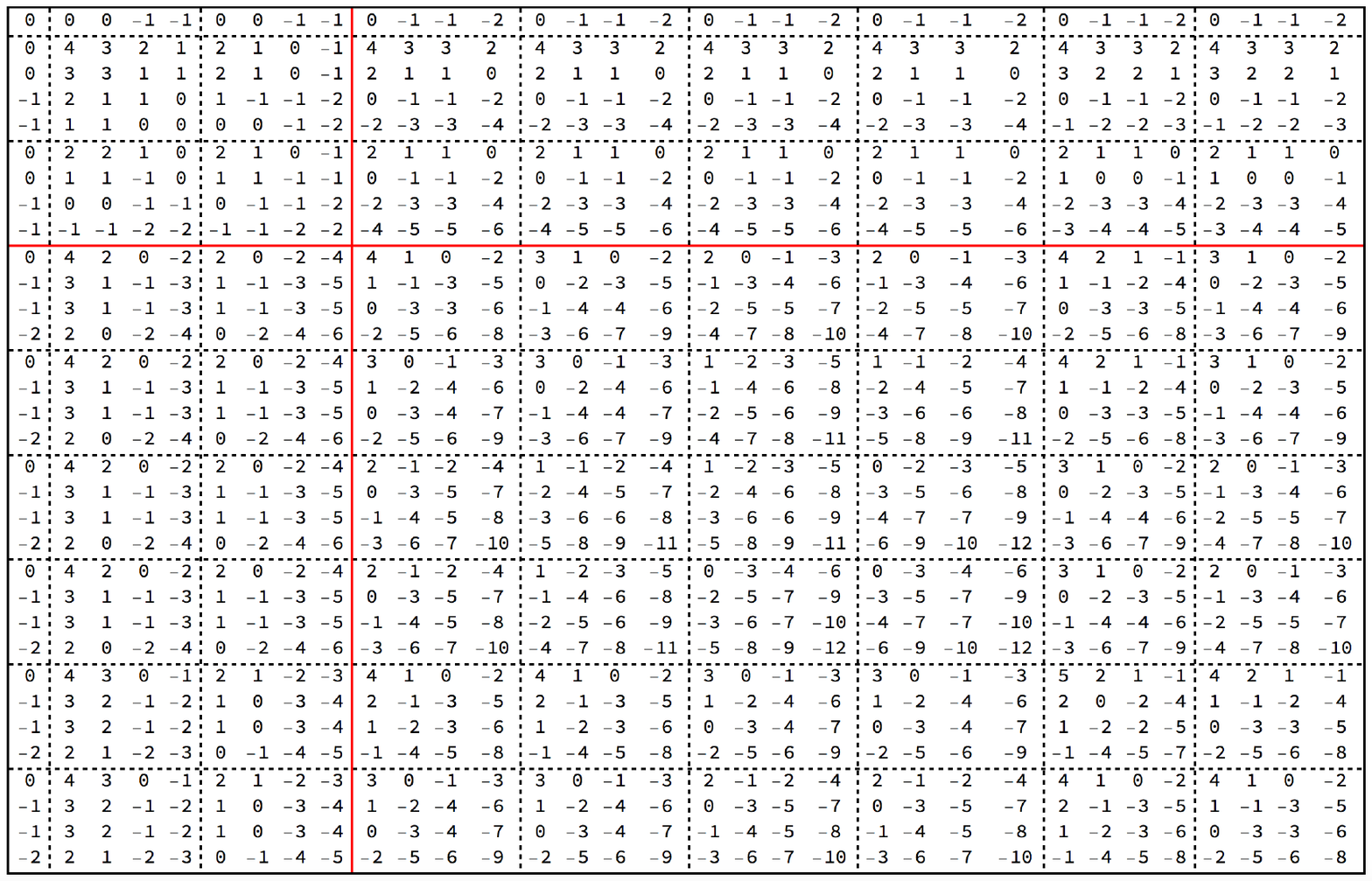}
  \caption{Quiver matrix for $\mathbf{9_{42}}$}  \label{942}
\end{figure}

The vectors corresponding to the linear terms in $a-$ and $q-$ degrees are (first 9 entries are for the level 1, and the next 24 for the level 2 generators, separated by semi-colon) :
$$\small{\!\!\!\!\!\!\!\!\!\!\!\!\mathbf{a}=(0,1,0,0,-1,1,0,0,-1;1,0,0,-1,0,-1,-1,-2,0,-1,-1,-2,-1,-2,-2,-3,1,0,0,-1,0,-1,-1,-2)},$$
$$\!\!\!\!\!\!\!\!\!\!\!\!\!\!\!\!\!\!\!\!\!\!\!\!\mathbf{q}=(0,-1,0,-2,-1,1,2,0,1;0,1,0,1,-2,-1,-2,-1,-3,-2,-3,-2,-5,-4,-5,-4,1,2,-1,0,-2,-1,-4,-3),$$
The matrix $C$ of the quiver is a symmetric 33 by 33 matrix. Note that it still has beautiful property that the diagonal entries are precisely the homological degrees (which are on the above formula for the colored polynomials seen as signs):
$$C_{i,i}=t_i,\quad i=1,\ldots,33.$$
Note that for level 1 nodes these are the corresponding homological gradings of the uncolored homology, while for the level 2 nodes, these are the homological gradings of the corresponding generators of the $S^2$-homology (they are precisely the $t_c$-gradings in the quadruply graded homology). More on this will be said in Section \ref{homology}.\\

To summarise, some of the properties of this quiver-form expression for $\mathbf{9_{42}}$, but also of all the quiver-forms for all knots from Section \ref{results}, are:\\

1.\quad  Correct values for $S^1$, $S^2$, $S^3$ and $S^4$ colored HOMFLY-PT polynomials.\footnote{In particular, we used \cite{NR} and \cite{paul}.}\\

2.\quad For $a=q^2$ gives the correct values for the colored Jones polynomial, for $r=1,\ldots,7$. (\cite{KnotAtlas}).\\

3.\quad It has the quiver form, perfectly working with the higher levels (higher exponents of $x$).\\

4.\quad Has the structure of "colored differentials" naturally built in -- the $q$-Pochhammers in the formula -- therefore suitable for categorification/refinement.

\subsection{Bottom row}
The values of the bottom rows of the colored HOMFLY-PT polynomials, that is coefficients of the lowest power of $a$ appearing in colored HOMFLY-PT polynomials, has also shown to be of independent interest. In particular bottom rows contain a significant non-trivial information about knots, including applications in combinatorics \cite{PSS}. The computations and all the previous results on bottom rows were relying on the fact that colored HOMFLY-PT polynomials for those knots have exponential growth. Since in this paper we have explicit formulas for colored HOMFLY-PT polynomials of knots with super-exponential growth, we can directly obtain the bottom rows for any color of those knots. Below we collect the values for the knot $\mathbf{9_{42}}$.

\subsubsection{Reduced colored HOMFLY-PT polynomial}
By (\ref{glpol}), the lowest $a$-degree of $P_r(\mathbf{9_{42}})(a,q)$ is $-3l-j+\alpha$. 
Together with the conditions $k+2l=r$, $0\le j\le k$ and $0\le\alpha\le l$, we have that the lowest degree is obtained for $\alpha=0$, $j=k$, and $l=[r/2]$.
The lowest $a$ degree is then $-r-l$.\\

Therefore, for even $r$, i.e. when $r=2l$, we have that the lowest $a$-degree is $-3l=-3r/2$, and it is equal to:
\[
r=2l:\quad P^{bottom}_r(\mathbf{9_{42}})(q)=q^{-6l^2+2l} \frac{(q;q)_{2l}}{(q;q)_l}=q^{-6l^2+2l} (q^{l+1};q)_l,
\]

while for odd $r$, i.e. when $r=2l+1$, we have that the lowest $a$-degree is $-3l-1=-3r/2-1/2$, and it is equal to:
\[
r=2l+1:\quad P^{bottom}_r(\mathbf{9_{42}})(q)= q^{-6l^2-6l} (q^{2l-1}+q^{4l+1})\frac{(q;q)_{2l+1}}{(q;q)_1(q;q)_l}.
\]

The table of lowest $a$-degrees for different $r$ is:
\begin{equation}
\begin{array}{c|c}
r & a_{\min}\\ \hline
1 & -1 \\
2& -3 \\
3 & -4\\
4 & -6\\
5 & -7\\
6 & -9 \\
7 & -10\\
\end{array}
\end{equation}

\subsubsection{Unreduced colored HOMFLY-PT polynomial}
For the unreduced version, we have to multiply the reduced polynomial by the unreduced polynomial for the unknot, which is:
\[
\bar{P}_r (U) (a,q)=a^{-r/2} q^{r/2} \frac{(a;q)_r}{(q;q)_r}.
\]
Since we are interested only in the bottom row here, this means that we just have to multiply by
$q^{r/2}/(q;q)_r$. Therefore:

For even $r$, $r=2l$, we have:
\[
r=2l:\quad \bar{P}^{bottom}_r(\mathbf{9_{42}})(q)= q^{-6l^2+2l} \frac{q^l}{(q;q)_l}.
\]

For odd $r$, $r=2l+1$, we have:
\[
r=2l+1:\quad \bar{P}^{bottom}_r(\mathbf{9_{42}})(q)= q^{-6l^2-6l} (q^{2l-1}+q^{4l+1})\frac{q^{l+1/2}}{(q;q)_1(q;q)_l}.
\]

\section{$S^r$-colored homology in (almost) quiver form of $\mathbf{9_{42}}$}\label{homology}

\subsection{Colored HOMFLY-PT homology of $\mathbf{9_{42}}$ -- first version}
A direct categorification/refinement of the above formulas, by using the analogy from the exponential growth cases (which have only level 1 nodes), would give quite reasonable guess for the homology. Those will have majority of the properties we predicted, even for the quadruply-graded case. Majority of the differentials would do to the job as predicted, in particular both canceling differentials, and all negative colored differentials. However, the other positive colored differentials would be slightly off. That can be fixed by adding certain number of generators that come in pairs, so that the added part of the superpolynomial would be of the form $(1+t) \times something$.

In particular for the $S^2$-colored homology, in $(a,q,t=t_c)$ grading, we would have:

\begin{eqnarray}
\!\!\!\mathcal{P}_2(\mathbf{9_{42}})(a,q,t)\!\!&\!\!=\!\!&\!\!1+ a (q^{-1} t^2+ q t^4) (1+q t^2) d_1 d_{-2}+ a^2 (t^8+ q^2 t^{10}+q^3 t^{12} +q^{6} t^{16})d_1 d_0 d_{-2} d_{-3}+\nonumber\\
&&\quad +\,(1+q  t) (1+q t^2) \left[a t^4 d_1 d_0 d_{-2} d_{-3} + a q t^5 d_1 d_{-2} d_{-3} \right]+\nonumber\\
&&\quad +\,(1+t) (1+q t^2)  a^2 q t^7  d_1 d_{-2} d_{-3},\label{hom942}
\end{eqnarray}
where 
\begin{eqnarray*}
d_1&=&1+a^{-1}q t^{-1},\\
d_0&=&1+a^{-1}t^{-3},\\
d_{-2}&=&1+a^{-1}q^{-2} t^{-5},\\
d_{-3}&=&1+a^{-1}q^{-3} t^{-7}.\\
\end{eqnarray*}
Note that the first two lines of $\mathcal{P}_2(\mathbf{9_{42}})(a,q,t)$ are natural refinement of the general formula (\ref{glpol}) for the colored HOMFLY-PT polynomial for $r=2$ case, with the first line corresponding to the (``square" of) the part coming from the level 1 generators/nodes, and the second line corresponds to the level 2 generators/nodes. Finally, the last line consists of those extra added terms that come in pairs (divisible by $1+t$). They behave like level 2 generators. 

Also note that even in the case of uncolored polynomial and homology of $\mathbf{9_{42}}$, the reduced HOMFLY-PT polynomial has seven terms, whereas for the homology has 9 terms (already Khovanov homology has 9 terms). The one extra pair of generators (exactly $1+t$ in the superpolynomial) was forced to be added in order to have two commuting differentials from the structural properties of the HOMFLY-PT homology, as predicted in \cite{DGR}.

\subsection{The first version of the $S^r$-colored homology of $\mathbf{9_{42}}$}
The explicit form of the quadruply-graded $S^r$-colored homology of $\mathbf{9_{42}}$, extending (\ref{hom942}) for every $r\in\N$, is given by (triply-graded homology is obtained by setting $Q=1$):
\begin{eqnarray}\label{ph942}
{\mathcal{\widetilde{P}}}_r(\mathbf{9_{42}})(a,q,t,Q)&=&\sum_{K+2L=r}\sum_{d_1+d_2+d_3=K} \sum_{d_4+d_5=L}\frac{(qt^2;qt^2)_r}{\prod_{i=1}^5 (qt^2;qt^2)_{d_i}}  a^{d_1+d_2+d_4+d_5} \times \\
&&\quad\times\quad (q Q)^{-d_1+d_2+d_5}  t^{2d_1+4d_2+4d_4+5d_5} (q t^2)^{-\frac{1}{2}(2d_1+4d_2+4d_4+5d_5)} (q t^2)^{\frac{1}{2}\sum_{i,j=1}^5 C_{ij}d_id_j} \times \nonumber\\
&&\quad\times\quad
(-a^{-1}q t^{-1} Q;(qt^2)^{-1})_{d_1+d_2+2d_4+d_5} (-a^{-1}q^{-r}t^{-2r-1}Q^{-1};(qt^2)^{-1})_{d_1+d_2+2d_4+2d_5}\times\nonumber\\
&&\quad\times\quad \frac{(-qt;qt^2)_{d_4}}{(qt^2;qt^2)_{d_4}} 
\frac{HP(q;d_1+d_2+2d_4)_{d_5}}{(qt^2;qt^2)_{d_5}}.\nonumber
\end{eqnarray}
Here $C$ is the symmetric matrix (\ref{mat942}) of the same quiver from the $S^r$-colored HOMFLY-PT polynomials expansion for $\mathbf{9_{42}}$
$$C=\left(
\begin{array}{ccc|cc}
2&2&0&2&2\\
2&4&0&4&4\\
0&0&0&0&0\\
\hline
2&4&0&4&4\\
2&4&0&4&5
\end{array}
\right)$$
and $HP(q;o)_m$ is the higher $a$-level categorification of the $q$-Pochhammer $(q;q)_m$.

\begin{definition}
The higher $a$-level categorification of the $q$-Pochhammer $(q;q)_m$ with shift $o$, denoted by $HP(q;o)_m$, is defined as:
\begin{equation}
HP(q;o)_m=\sum_{i=0}^m a^i q^{-i} Q^{-i} {m\brack i}_{qt^2}
(qt^2)^{(m+o)i} (-qt;qt^2)_{m-i} (-t;qt^2)_i.
\end{equation}
\end{definition}

\subsubsection{Properties of $HP(q;o)_m$}
Triply-graded version is obtained by setting $Q=1$. We clearly have that $HP(q;o)_m$ categorifies $(q;q)_m$:
\[
HP(q;o)_m\mid_{{}_{t=-1}}=(q;q)_m.
\]
Moreover, it has very interesting behaviour with respect to the homology of differential $d_n$. To that end, we set:
$$D_i=1+a q^{-i} t^{3-2i} Q^{-1}, \quad i\in \Z.$$
\begin{definition}
For three polynomials $P_1,P_2, D\in \N[a^{\pm},q^{\pm},t^{\pm},Q^{\pm}]$, we say that 
\[
P_1\mod D = P_2,
\]
if 
\[P_1(a,q,t,Q)-P_2(a,q,t,Q) = D(a,q,t,Q) Y(a,q,t,Q),\] 
for some polynomial $Y\in\N[a^{\pm},q^{\pm},t^{\pm},Q^{\pm}]$ with non-negative integer coefficients.
\end{definition}
Then we have:
\begin{lemma}
\begin{eqnarray}
&&HP(q;o)_m \mod D_{2-m-o} = \prod_{\ell=1}^m D_{2-\ell-m-o}.\\
&&HP(q;o)_m \mod D_{1-m-o} = q^{m(m+1)/2} t^{m^2} \prod_{\ell=1}^m D_{1+\ell-m-o}.
\end{eqnarray}
\end{lemma}

\subsubsection{Properties of the $S^r$-colored homology }

\paragraph{\underline{Positive colored differentials}}
\begin{proposition}
For every $l=0,\ldots,r-1$, we have:
\[
{{\mathcal{\widetilde{P}}}}_r(\mathbf{9_{42}})(a,q,t,Q)-{\mathcal{\widetilde{P}}}_l(\mathbf{9_{42}})(a,q,t,Q (qt^2)^{r-l})=(1+a^{-1}q^{1-l} t^{-1-2l} Q)  X_+(a,q,t,Q),
\]
where $X_+\in \N[a^{\pm},q^{\pm},t^{\pm},Q^{\pm}]$ is a polynomial with non-negative integer coefficients.
\end{proposition}
With the notation from above this can be written as:
$${\P}_r(\mathbf{9_{42}})(a,q,t,Q) \mod D_{1-l}   =  {\P}_l(\mathbf{9_{42}})(a,q,t,Q (qt^2)^{r-l}),  \quad\quad l=0,\ldots,r-1,$$
exactly as predicted in [GGS].\\

In terms of $(a,q,t_r,t_c)$-grading (note that $t_c=t$, and $t_r=1$) gives triply-graded homology, we have:
$${\P '}_r(\mathbf{9_{42}})(a,q,t_r,t_c)\mod\left(1+a q^{l-1} t_r t_c^{2l+1}\right)   =  {\P '}_l(\mathbf{9_{42}})\left(a,q (qt_c^2)^{r-l},t_r (qt_c^2)^{\frac{r-l}{2}},t_c (qt_c^2)^{-\frac{r-l}{2}}\right),$$ for $ l=0,\ldots,r-1.  $

\paragraph{\underline{Negative colored differentials}}
\begin{proposition}
For every $l=0,\ldots,r-1$, we have:
\[
{{\mathcal{\widetilde{P}}}}_r(\mathbf{9_{42}})(a,q,t,Q)-{\mathcal{\widetilde{P}}}_l(\mathbf{9_{42}})(a,q,t,Q)=(1+a^{-1}q^{-r-l} t^{-1-2r-2l} Q^{-1})  X_-(a,q,t,Q),
\]
where $X_-\in \N[a^{\pm},q^{\pm},t^{\pm},Q^{\pm}]$ is a polynomial with non-negative integer coefficients.
\end{proposition}
With the notation from above this can be written as:
$${\P}_r(\mathbf{9_{42}})(a,q,t,Q) \mod (1+a q^{r+l} t^{2r+2l+1}Q)   =  {\P}_l(\mathbf{9_{42}})(a,q,t,Q), \quad\quad l=0,\ldots,r-1,$$
exactly as predicted in [GGS]. Also, note that this statement can be written explicitely in the three-grading setting, by putting simply $Q=1$.\\

In terms of $(a,q,t_r,t_c)$-grading (note that $t_c=t$, and $t_r=1$ gives triply-graded homology), we have:
$${\P}'_r(\mathbf{9_{42}})(a,q,t_r,t_c) \mod \left(1+a q^{r+l} t_r^3 t_c^{2r+2l+1}\right)   =  {\P}'_l(\mathbf{9_{42}})(a,q ,t_r,t_c) \quad\quad l=0,\ldots,r-1.$$

\paragraph{\underline{Self-symmetry}}
Self-symmetry is most elegantly stated in $(a,Q,t_r,t_c)$ four-grading:
\begin{proposition}
For every $r\in\N$:
\[
{\P}''_r(\mathbf{9_{42}})(a,Q,t_r,t_c)={\P}''_r(\mathbf{9_{42}})(a,Q^{-1} t_r^{-2} t_c^{-2r} ,t_r,t_c)
\]
\end{proposition}

\subsubsection{Size of the homology, and of the quiver}
From the expansion (\ref{ph942}) of $\P_r(\mathbf{9_{42}})$, one can obtain the size of the reduced $S^r$-colored HOMFLY-PT homology of $\mathbf{9_{42}}$, for any $r\in\N$:
\begin{eqnarray*}
\dim \mathcal{H}_r(\mathbf{9_{42}})&=&\P_r(\mathbf{9_{42}})(a=1,q=1,t=1,Q=1)=\\
&=&\quad\sum_{K+2L=r}\sum_{d_1+d_2+d_3=K} \sum_{d_4+d_5=L}\binom{r}{K}\binom{2L}{L}\binom{K}{d_1,d_2,d_3}\binom{L}{d_4,d_5}^2   \times\\
&&\quad\quad\quad\quad\quad\quad\quad\quad\quad\quad\quad\quad\quad \times\,\, 2^{d_1+d_2+2d_4+d_5} 2^{d_1+d_2+2d_4+2d_5} 2^{d_4} 4^{d_5}=\\
&=&\quad \sum_{K+2L=r}\sum_{d_1+d_2+d_3=K} \sum_{d_4+d_5=L}\binom{r}{K}\binom{2L}{L}\binom{K}{d_1,d_2,d_3}\binom{L}{d_4,d_5}^2 4^{d_1+d_2} 32^{d_4+d_5}=\\
&=&\quad \sum_{K+2L=r}\binom{r}{K}\binom{2L}{L} 32^L \left(\sum_{d_4=0}^L \binom{L}{d_4}^2\right)\sum_{d_3=0}^K\binom{K}{d_3}\sum_{d_1+d_2=K-d_3}\binom{K-d_3}{d_1,d_2} 4^{d_1+d_2}=\\
&=&\quad \sum_{K+2L=r}\binom{r}{K}\binom{2L}{L} 32^L \binom{2L}{L}\sum_{d_3=0}^K\binom{K}{d_3} 4^{K-d_3} 2^{K-d_3}=\\
&=&\quad \sum_{K+2L=r}\binom{2L}{L}^2 32^L \sum_{d_3=0}^K\binom{K}{d_3} 8^{d_3} =\, \sum_{K+2L=r}\binom{r}{K}\binom{2L}{L}^2 32^L \sum_{d_3=0}^K\binom{K}{d_3} 8^{d_3}=\\
&=&\quad \sum_{K+2L=r}\binom{r}{K}\binom{2L}{L}^2 \,9^K \,32^L.
\end{eqnarray*}
The sizes of the first couple of colored homologies of $\mathbf{9_{42}}$ are therefore:
$$9,\,209,\,4185,\,105633,\,2651049,\, 71025521,\ldots$$

\vskip 0.4cm

As for the quiver, for the refined/homological version, we actually have almost quiver form. In particular, this is due to the presence of the categorification of $q$-Pochhammer symbols, which are done in two different ways, corresponding to the summation variables $d_4$ and $d_5$, in the last line of (\ref{ph942}). Nevertheless, we still have the generators only of levels ($x$-degrees) 1 and 2. 
There are 9 generators at level 1 (each summation variables $d_1$ and $d_2$ carries 4 generators, and a single generator corresponding to variable $d_3$), and 64 generators at level 2 (or 32 or 128, depending how one counts).

\subsection{The alternative version of $S^2$ homology of $\mathbf{9_{42}}$}

This $S^2$ homology is slightly different comparing to the first version -- (\ref{hom942}), and its main advantage that it an addition to other properties, it also gives the correct $SO(3)$ homology. The only property it "looses" is the universal differential (which might be the right thing to do with higher levels generators). The "essential" part of the homology is the same, however added pairs of generators, that is, part of the form $(1+t)\times something$, is different than in (\ref{hom942}).
\begin{eqnarray} 
\mathcal{P}^{S^2}_{red}(\mathbf{9_{42}})&=&1+a(q^{-1}t^2+qt^4)(1+q^2t)d_1d_{-2}+a^2(t^8+q^2t^{10}+q^3t^{12}+q^6t^{16})d_1d_{0}d_{-2}d_{-3}+\nonumber\\
&&+\,(1+q^2t^3)(at^4d_1d_0d_{-2}d_{-3}+aqt^5d_1d_{-2}d_{-3})+ \label{formula}\\
&&+\,(1+t)a^2qt^7d_1d_{-2}d_{-3}+(1+t)(a^3q^{-1}t^8+a^3 q^5 t^{16} + a^4 q^3 t^{17}) d_1d_0d_{-2}d_{-3}
\nonumber
\end{eqnarray}
Here again
\begin{eqnarray*}
d_1&=& 1+a^{-1}qt^{-1}\\
d_0&=& 1+a^{-1}t^{-3}\\
d_{-2}&=& 1+a^{-1}q^{-2}t^{-5}\\
d_{-3}&=& 1+a^{-1}q^{-3}t^{-7}\end{eqnarray*}

One can check that indeed this superpolynomial satisfies  the following conditions:

\begin{itemize}
\item{\quad { it gives the correct colored polynomial for $t=-1$}}
\item{\quad differentials $d_1$ and $d_{-2}$ are canceling}
\item{\quad the homology with respect to $d_{-3}$ is isomorphic to uncolored homology, without any grading change in $a$ and $q$}
\item{\quad the homology with respect to $d_{0}$ is also isomorphic to uncolored homology, precisely with the predicted shift in $q$ grading: $q\to q^2$}
\item{\quad one can also add the fourth grading $t_r$ or $Q$ ($t$ in the formula above is $t_c$), so that all of the properties above together with the self symmetry hold on quadruply graded level exactly as predicted in [GGS].}\footnote{The quadruply graded version is:
\begin{eqnarray*} 
\mathcal{\tilde{P}}^{S^2}_{red}(\mathbf{9_{42}})(a,q,t_r,t_c=t)\!\!\!&\!\!\!\!\!\!=\!\!\!\!\!\!&\!\!\!1+a(q^{-1}t_r^2 t^2+t_r^2 t^4+q t_r^4 t^4+q^2t_r^4 t^6)D_1D_{-2}+a^2(t_r^4 t^8+q^2 t_r^6 t^{10}+q^3t_r^6 t^{12}+q^6 t_r^8 t^{16})D_1D_{0}D_{-2}D_{-3}\\\nonumber
\!\!\!&\!\!\!\!\!\!&\!\!\!+ (1+q^2t_r^{-1}t^3)(at_r^4 t^4D_1D_0D_{-2}D_{-3}+aq t_r^5 t^5 D_1D_{-2}D_{-3})+(1+t_r t)a^2q t_r^5 t^7D_1D_{-2}D_{-3}\\ \nonumber
\!\!\!&\!\!\!\!\!\!&\!\!\!+ (1+t_r t)(a^3q^{-1}t_r^6 t^8+a^3 q^5 t_r^{10} t^{16} + a^4 q^3 t_r^{11} t^{17}) D_1D_0D_{-2}D_{-3}\label{formula4}\end{eqnarray*}

Here \begin{eqnarray*}
D_1&=& 1+a^{-1}q t_r^{-1}t^{-1}\\
D_0&=& 1+a^{-1}t_r^{-1}t^{-3}\\
D_{-2}&=& 1+a^{-1}q^{-2}t_r^{-3}t^{-5}\\
D_{-3}&=& 1+a^{-1}q^{-3}t_r^{-3}t^{-7}.\end{eqnarray*}
In particular, if we set 
$$\tilde{P}_{red}^{S^2}(\mathbf{9_{42}})(a,Q,t_r,t_c):= \mathcal{\tilde{P}}^{S^2}_{red}(\mathbf{9_{42}})(a,q,t_r q^{-1/2},t_c q^{1/2})_{\mid q \to Q},$$
then self-symmetry of the quadruply-graded homology can be written as:
$$\tilde{P}_{red}^{S^2}(\mathbf{9_{42}})(a,Q,t_r,t_c)=\tilde{P}_{red}^{S^2}(\mathbf{9_{42}})(a,Q^{-1},t_r Q^{-1},t_c Q^{-2}).$$}

\item{\quad It gives the correct $SO(3)$ homology that matches the results of Cooper, Hogancamp and Krushkal (\cite{CHK,ben}), as we explicitely demonstrate below in Section \ref{so3}.}

\item{\quad Finally, there exists a differential $d_{-4}$ of $(a,q,t_r,t)$ four-degree $(-1,-4,-5,-9)$, so that the the homology of $\mathcal{H}^{S^2}(\mathbf{9_{42}})$ with respect to $d_{-4}$ is isomorphic to the SO(6) homology, which in turn can be obtained from the triply-graded Kauffman homology with respect to the SO(6) differential ($d_5$ in that notation). This lifts the property that the $(sl(4),\Lambda^2)$ representation is isomorphic to the vector representation of $so(6)$
(we implicitly use the property that $S^2$ and $\Lambda^2$ homologies are mirror symmetric with respect to $q\leftrightarrow q^{-1}$). The Kauffman homology and its properties in general will be pursued in more details in the future paper.}
\end{itemize}

\subsubsection{$SO(3)$ homology of $\mathbf{9_{42}}$}\label{so3}

We shall obtain $SO(3)$ homology of $\mathbf{9_{42}}$ as the $(sl(2),S^2)$ homology, i.e. as the homology of the $S^2$ unreduced homology with respect to the differential $d_2$, of $(a,q,t)$ tri-grading $(-1,2,-1)$. The unreduced homology is obtained as the tensor product of the reduced homology of a knot and the unreduced homology of the unknot whose super polynomial is given by 
$${\mathcal{\bar{P}}}^{S^2}(U)=a^{-1}q\frac{(1+at)(1+aqt^3)}{(1-q)(1-q^2t^2)}.$$
In addition, in order to match the computations of Cooper, Hogancamp and Krushkal, the knot they observe is the mirror image of $\mathbf{9_{42}}$ for which we obtained the formula (\ref{formula}), therefore the super polynomial of the unreduced homology of $\overline{\mathbf{9_{42}}}$ is just 
\begin{equation}\label{miror}
\mathcal{\bar{P}}^{S^2}(\overline{\mathbf{{9}_{42}}})(a,q,t)=\mathcal{{P}}^{S^2}_{red}(\mathbf{9_{42}})(a^{-1},q^{-1},t^{-1})\mathcal{\bar{P}}^{S^2}(U).
\end{equation}
Finally, we have to compute the homology with respect to $d_2$, i.e. the value of this super polynomial modulo $d_2=1+a^{-1}q^2t^{-1}$. For the unknot part we have:
$$\mathcal{\bar{P}}^{S^2}(U)=a^{-1}q (1+q +\frac{q^2t^2+a q^2t^3}{1-q^2t^2}) \quad (\,\mathop{\rm mod}\, d_2).$$
In addition, for each pair of generators of reduced homology that are cancelled by $d_{-2}$ differential, the corresponding unreduced polynomial modulo $1+a^{-1}q^2t^{-1}$ becomes finite. 
More precisely:
$$d_{-2} \mathcal{\bar{P}}^{S^2}(U) =(1+a^{-1}q^{-2}t^{-5}) a^{-1}q (1+q +\frac{q^2t^2+a q^2t^3}{1-q^2t^2}) = a^{-2}q^{-1}t^{-5} (1+q+q^2t^2+aq^2t^3+aq^3t^5+aq^4t^5)  \quad (\,\mathop{\rm mod}\,  d_2).$$
Moreover:
$$d_1d_{-2} \mathcal{\bar{P}}^{S^2}(U) = a^{-3}t^{-6} (1+q^2t^2)(1+at)(1+aqt^3)  \quad (\,\mathop{\rm mod}\,  d_2),$$
and
$$d_1d_0 d_{-2} \mathcal{\bar{P}}^{S^2}(U) = a^{-4}t^{-9} (1+at)(1+aqt^3)(1+aq^2t^5)  \quad (\,\mathop{\rm mod}\,  d_2),$$
as well as
$$d_1d_0 d_{-2}d_{-3} \mathcal{\bar{P}}^{S^2}(U) = a^{-5}q^{-3}t^{-16} (1+at)(1+aqt^3)(1+aq^2t^5)(1+aq^3t^7)  \quad (\,\mathop{\rm mod}\,  d_2).$$
Finally, from (\ref{formula}) we have (equalities are $\,\mathop{\rm mod}\,  d_2$):
\begin{eqnarray*}
\mathcal{P}^{S^2}_{red}(\mathbf{9_{42}})&=&1+a(q^{-1}t^2+t^4+qt^4+q^2t^6)d_1d_{-2}+a^2(q^2t^{10}+q^3t^{12}+q^6t^{16})d_1d_{0}d_{-2}d_{-3}+\\\nonumber
&&+at^4d_1d_0d_{-2}d_{-3}+aqt^5d_1d_{-2}d_{-3}+aq^3t^8d_1d_{-2}d_{-3}+(a^2qt^7+a^2q t^8) d_1d_{-2}d_{-3}+\\\nonumber
&&+(a^3q^{-1}t^8+a^3q^{-1}t^9)(1+a^{-1}t^{-3})d_1d_{-2}d_{-3}+(1+t)a^4 q^3 t^{17}(1+a^{-1}q^2t^{-1}) d_1d_0d_{-2}d_{-3}+\\
&&+a^2t^{8}d_1d_{0}d_{-2}d_{-3}+a q^2 t^7d_1d_0d_{-2}d_{-3}=\\
&=&1+a(q^{-1}t^2+t^4+qt^4+q^2t^6)d_1d_{-2}+a^2(q^2t^{10}+q^3t^{12}+q^6t^{16})d_1d_{0}d_{-2}d_{-3}+\\\nonumber
&&+at^4d_1d_0d_{-2}d_{-3}+(aqt^5+aq^3t^8+a^2qt^7+a^2q t^8)d_1d_{-2}d_{-3}+\\\nonumber
&&+(a^3q^{-1}t^8+a^3q^{-1}t^9+a^2q^{-1}t^5+a^2q^{-1}t^{6})d_1d_{-2}d_{-3}+(1+t)a^4 q^3 t^{17}(1+a^{-1}q^2t^{-1}) d_1d_0d_{-2}d_{-3}+\\
&&+a^2 t^8 (1+a^{-1}q^2t^{-1})d_1d_{0}d_{-2}d_{-3}=\\
&=&1+a(q^{-1}t^2+t^4+qt^4+q^2t^6)d_1d_{-2}+a^2(q^2t^{10}+q^3t^{12}+q^6t^{16})d_1d_{0}d_{-2}d_{-3}+\\\nonumber
&&+(aq^3t^8)d_1d_{-2}d_{-3}+(a^2q^{-1}t^5)d_1d_{-2}d_{-3}=\\
&=&1+a(q^{-1}t^2+t^4+qt^4+q^2t^6)d_1d_{-2}+(a^2q^2t^{10}+a^2q^3t^{12}+a^2q^6t^{16}+at^4)d_1d_{0}d_{-2}d_{-3}+\\
&&+(aq^3t^8+t^1)d_1d_{-2}+(a^2q^{-1}t^5+a q^{-4} t^{-2})d_1d_{-2}=\\
&=&1+(a^2q^2t^{10}+a^2q^3t^{12}+a^2q^6t^{16}+at^4)d_1d_{0}d_{-2}d_{-3}+\\
&&+(aq^{-1}t^2+at^4+aqt^4+aq^2t^6+aq^3t^8+t^1+a^2q^{-1}t^5+a q^{-4} t^{-2})d_1d_{-2}=\\
&=&1+(a^2q^2t^{10}+a^2q^3t^{12}+a^2q^6t^{16}+at^4)d_1d_{0}d_{-2}d_{-3}+\\
&&+(aq^{-1}t^2+at^4+aq^2t^6+aq^3t^8+t^1+a q^{-4} t^{-2})d_1d_{-2} \quad (\,\mathop{\rm mod}\,  d_2).
\end{eqnarray*}
By returning this into (\ref{miror}) and by using the expressions of $d_1d_{-2}$ and $d_1d_{0}d_{-2}d_{-3}$ modulo $d_2$, we obtain:
\begin{eqnarray*}
&\mathcal{\bar{P}}^{S^2}(\overline{\mathbf{9_{42}}})(a,q,t)=a^{-1}q (1+q+q^2t^2+a q^2t^3 + \frac{q^4t^4+aq^4t^5}{1-q^4t^2})+\\
&\!\!\!\!\!\!+a^{-1}q (a^{-2}q^{-6}t^{-16}+a^{-2}q^{-3}t^{-12}+a^{-2}q^{-2}t^{-10}+a^{-1}t^{-4})(1+at)(1+aqt^3)(1+aq^2t^5)(1+aq^3t^7)+\\
&\!\!\!\!\!\!+a^{-1}q (a^{-1}qt^{-2}+a^{-1}t^{-4}+a^{-1}q^{-2}t^{-6}+a^{-1}q^{-3}t^{-8}+t^{-1}+a^{-1} q^{4} t^{2})(1+q^2t^2)(1+at)(1+aqt^3).
\end{eqnarray*}
Finally, since 
\begin{eqnarray*}
(a^{-1}t^{-4}+t^{-1})(1+q^2t^2)=t^{-1}+q^2 t+a^{-1}t^{-4}+a^{-1}q^2t^{-2}=q^2 t +a^{-1}t^{-4}\quad (\,\mathop{\rm mod}\,  d_2),
\end{eqnarray*}
we have that, modulo $d_2$:
\begin{eqnarray*}&&\mathcal{\bar{P}}^{S^2}(\overline{\mathbf{9_{42}}})(a,q,t)=a^{-1}q (1+q+q^2t^2+a q^2t^3 + \frac{q^4t^4+aq^4t^5}{1-q^4t^2})+\\
&&\!\!\!\!\!\!\!\!\!\!\!\!+a^{-1}q (a^{-2}q^{-6}t^{-16}+a^{-2}q^{-3}t^{-12}+a^{-2}q^{-2}t^{-10}+a^{-1}t^{-4})(1+at)(1+aqt^3)(1+aq^2t^5)(1+aq^3t^7)+\\
&&\!\!\!\!\!\!+a^{-1}q (a^{-1}qt^{-2}+a^{-1}q^3+a^{-1}q^{-2}t^{-6}+a^{-1}t^{-4}+a^{-1}q^{-3}t^{-8}+a^{-1}q^{-1}t^{-6})(1+at)(1+aqt^3)+\\
&&\!\!\!\!\!\!+(a^{-1} q^{4} t^{2}+a^{-1}q^6 t^4+q^2t+a^{-1}t^{-4})(1+at)(1+aqt^3).
\end{eqnarray*}

Finally, after collapsing to bi-grading by setting $a=q^2$, we obtain the superpolynomial of the unreduced $(sl(2),S^2)$-homology:
\begin{eqnarray*}\mathcal{\bar{P}}^{sl(2),S^2}(\overline{\mathbf{9_{42}}})(q,t)&=& (q^{-1}+1+q t^2+q^3t^3 + \frac{q^3t^4+q^5t^5}{1-q^4t^2})+\\
&&\!\!\!\!\!\!\!\!\!\!\!\!+(q^{-11}t^{-16}+q^{-8}t^{-12}+q^{-7}t^{-10}+q^{-3}t^{-4})(1+q^2t)(1+q^3t^3)(1+q^4t^5)(1+q^5t^7)+\\
&&\!\!\!\!\!\!+(q^{-2}t^{-2}+1+q^{-5}t^{-6}+q^{-3}t^{-4}+q^{-6}t^{-8}+q^{-4}t^{-6})(1+q^2t)(1+q^3t^3)+\\
&&\!\!\!\!\!\!+( q t^{2}+q^3 t^4+q t+q^{-3}t^{-4})(1+q^2t)(1+q^3t^3).
\end{eqnarray*}

The last homology, up to rescaling $q \to q^2$, and an overall shift, matches the result for the $SO(3)$ homology (defined in \cite{CHK}) of $\mathbf{9_{42}}$, computed by their computer program \cite{ben}.\footnote{
I would like to thank Ben Cooper for sharing this computation. For $\mathbf{9_{42}}$, the Poincare polynomial of the $SO(3)$ homology is given by:
\begin{eqnarray*}
\scriptstyle{\mathcal{P}^{SO(3)}(\mathbf{9_{42}})(q,t)}&\scriptstyle{=}&\scriptstyle{1 + 2/q^2 + 3 q^2 + 1/(q^{26} t^{18}) + 1/(q^{22} t^{17}) + 1/(q^{20} t^{15}) + 1/(q^{20} t^{14}) + 1/(q^{16} t^{14}) + 1/(q^{18} t^{13}) + 1/(q^{16} t^{13})+ }  \\
 &&\scriptstyle{+1/(q^{18} t^{12}) + 1/(q^{14} t^{12}) + 1/(q^{16} t^{11}) + 2/(q^{14} t^{11}) + 1/(q^{16} t^{10}) + 2/(q^{12} t^{10}) + 1/(q^{10} t^{10}) + 3/(q^{12} t^9) + 1/(q^8 t^9) + }\\
 &&\scriptstyle{+1/(q^{14} t^8) + 1/(q^{12} t^8) + 1/(q^{10} t^8) + 2/(q^8 t^8) + 4/(q^{10} t^7) + 1/(q^8 t^7) + 1/(q^6 t^7) + 3/(q^{10} t^6) +} \\
 &&\scriptstyle{+1/(q^8 t^6) + 4/(q^6 t^6) + 2/(q^8 t^5) + 4/(q^6 t^5) + 1/(q^4 t^5) + 1/(q^2 t^5) + 1/(q^8 t^4) + 4/(q^4 t^4) + 1/(q^2 t^4) + 2/t^3 + 4/(q^4 t^3) +}\\
&&\scriptstyle{+ 1/(q^2 t^3) + 3/t^2 + 1/(q^6 t^2) + 2/(q^4 t^2) + 2/(q^2 t^2) + q^2/t^2 + 1/t + 3/(q^2 t) + (2 q^2)/t + t + 3 q^2 t + 2 q^4 t +} \\
&&\scriptstyle{+q^2 t^2 + 3 q^4 t^2 + q^6 t^2 + q^8 t^2 + 
 q^4 t^3 + 2 q^6 t^3 + 2 q^8 t^3 + q^6 t^4 + q^8 t^4 +} \\
&&\scriptstyle{+q^{10} t^4 + q^8 t^5 + q^{10} t^5 + q^8 t^6 + q^{12} t^6 + q^{12} t^7 + q^{14} t^9 + q^{18} t^{10}+(t^2 q^2 + t^3 q^6)/(1-t^2 q^4).}
 \end{eqnarray*}
 }


\appendix

\section{Non-uniqueness of a quiver for knot $\mathbf{5_1}$}\label{ap}
For knot $\mathbf{5_1}$, the data for the ansatz quiver-like form (\ref{specred}) is given in (\ref{mat51}), i.e. we have
\begin{eqnarray}
P_r({\bf 5_{1}})(a,q)&=&\sum_{d_1+d_2+d_3=r} \frac{(q;q)_r}{(q;q)_{d_1}(q;q)_{d_2}(q;q)_{d_3}} (-1)^{d_2+d_3} a^{2d_1+3d_2+3d_3} q^{-2d_1-d_2+d_3} q^{-\frac{1}{2}(3d_2+5d_3)}\times \nonumber \\
&&\quad\times\quad  q^{\frac{1}{2}(3d_2^2+5d_3^2+2d_1d_2+6d_1d_3+8d_2d_3)} \,\,
(a^{-1}q;q^{-1})_{d_2+d_3} .\label{for51}
\end{eqnarray}
In order to get the quiver-form, we need to rewrite this expression in order not to have the last $q$-Pochhammer $(a^{-1}q;q^{-1})_{d_2+d_3}$. As explained in Section \ref{nonun}, we can do it in two ways:
\begin{equation}\label{op1}
(a^{-1}q;q^{-1})_{d_2+d_3}=(a^{-1}q;q^{-1})_{d_2} (a^{-1}q^{1-d_2};q^{-1})_{d_3},
\end{equation}
or
\begin{equation}\label{op2}
(a^{-1}q;q^{-1})_{d_2+d_3}=(a^{-1}q;q^{-1})_{d_3} (a^{-1}q^{1-d_3};q^{-1})_{d_2}.
\end{equation}
If we use (\ref{op1}), by $q$-binomial identity, we have:
\begin{eqnarray}
(a^{-1}q;q^{-1})_{d_2+d_3}&=&(a^{-1}q^{2-d_2};q)_{d_2} (a^{-1}q^{2-d_2-d_3};q)_{d_3}=\nonumber\\
&=&\sum_{e_2+e_3=d_2}(-1)^{e_2} a^{-e_2} q^{(2-e_2-e_3)e_2}q^{\frac{1}{2}(e_2^2-e_2)}\frac{(q;q)_{d_2}}{(q;q)_{e_2}(q;q)_{e_3}}\times\nonumber\\
&&\quad\times \sum_{e_4+e_5=d_3}(-1)^{e_4} a^{-e_4} q^{(2-e_2-e_3-e_4-e_5)e_4}q^{\frac{1}{2}(e_4^2-e_4)}\frac{(q;q)_{d_3}}{(q;q)_{e_4}(q;q)_{e_5}}.\label{ver1}
\end{eqnarray}
On the other hand if we use (\ref{op2}), and $q$-binomial identity, with the same change of variables, $d_2=e_2+e_3$ and $d_4=e_4+e_5$, we would obtain:
\begin{eqnarray}
(a^{-1}q;q^{-1})_{d_2+d_3}&=&(a^{-1}q^{2-d_2-d_3};q)_{d_2} (a^{-1}q^{2-d_3};q)_{d_3}=\nonumber\\
&=&\sum_{e_2+e_3=d_2}(-1)^{e_2} a^{-e_2} q^{(2-e_2-e_3-e_4-e_5)e_2}q^{\frac{1}{2}(e_2^2-e_2)}\frac{(q;q)_{d_2}}{(q;q)_{e_2}(q;q)_{e_3}}\times\nonumber\\
&&\quad\times \sum_{e_4+e_5=d_3}(-1)^{e_4} a^{-e_4} q^{(2-e_4-e_5)e_4}q^{\frac{1}{2}(e_4^2-e_4)}\frac{(q;q)_{d_3}}{(q;q)_{e_4}(q;q)_{e_5}}.\label{ver2}
\end{eqnarray}
By replacing (\ref{ver1}) into (\ref{for51}), with $e_1=d_1$, we get
\begin{eqnarray}
P_r({\bf 5_{1}})(a,q)&=&\sum_{e_1+e_2+e_3+e_4+e_5=r} \frac{(q;q)_r}{\prod_{i=1}^5(q;q)_{e_i}}  (-1)^{e_3+e_5}a^{2e_1+2e_2+3e_3+2e_4+3e_5} q^{-2e_1-e_3+e_5+2e_4} \times \nonumber \\
&&\quad\quad\times\quad q^{-\frac{1}{2}(2e_2+3e_3+4e_4+5e_5)}  q^{\frac{1}{2}\sum_{i,j}C^{(1)}_{i,j}e_ie_j},\label{for511}
\end{eqnarray}
where
\begin{equation}\label{C151}
C^{(1)}=\left(\begin{array}{ccccc}
0&1&1&3&3\\
1&2&2&3&4\\
1&2&3&3&4\\
3&3&3&4&4\\
3&4&4&4&5
\end{array}
\right).
\end{equation}
On the other hand, by replacing (\ref{ver2}) into (\ref{for51}), again with $e_1=d_1$, we get
\begin{eqnarray}
P_r({\bf 5_{1}})(a,q)&=&\sum_{e_1+e_2+e_3+e_4+e_5=r} \frac{(q;q)_r}{\prod_{i=1}^5(q;q)_{e_i}}   (-1)^{e_3+e_5} a^{2e_1+2e_2+3e_3+2e_4+3e_5} q^{-2e_1-e_3+e_5+2e_4} \times \nonumber \\
&&\quad\quad\times\quad  q^{-\frac{1}{2}(2e_2+3e_3+4e_4+5e_5)} q^{\frac{1}{2}\sum_{i,j}C^{(2)}_{i,j}e_ie_j},\label{for512}
\end{eqnarray}
where
\begin{equation}\label{C251}
C^{(2)}=\left(\begin{array}{ccccc}
0&1&1&3&3\\
1&2&2&3&3\\
1&2&3&4&4\\
3&3&4&4&4\\
3&3&4&4&5
\end{array}
\right).
\end{equation}
In this way, we obtain two different quiver matrices $C^{(1)}$ and $C^{(2)}$ corresponding 
to the knot $\mathbf{5_1}$, directly from the ansatz form (\ref{for51}).  
We note that both of them are part of the permutahedron structure of the quivers with minimal size corresponding to the knot $\mathbf{5_1}$ obtained in \cite{permutahedron}.




\end{document}